\documentclass[12pt]{article}
\title{Poisson deformations of affine symplectic 
varieties}
\author{Yoshinori Namikawa}
\date{ }
\usepackage{amsmath,amssymb,amsthm, amscd}
\chardef\bslash=`\\

\def\0{{\mathcal O}}
\def\g{{\mathfrak g}}
\def\p{{\mathfrak p}}
\def\h{{\mathfrak h}}

\def\l{{\mathfrak l}}
\def\n{{\mathfrak n}}
\def\b{{\mathfrak b}}

\begin{document}
\maketitle

\begin{center}
{\bf Introduction} 
\end{center}
\vspace{0.2cm}
 
A symplectic variety $X$ is a normal algebraic variety 
(defined over $\mathbf{C}$) which admits an everywhere 
non-degenerate d-closed 2-form  
$\omega$ on the regular locus $X_{reg}$ of 
$X$ such that, for any resolution $f: \tilde{X} \to X$ 
with $f^{-1}(X_{reg}) \cong X_{reg}$, the 2-form $\omega$ extends to a regular 
closed 2-form on $\tilde{X}$ (cf. [Be]). There is a natural Poisson structure 
$\{\; , \;\}$ on $X$ determined by $\omega$. 
Then we can introduce the notion of a Poisson deformation 
of $(X, \{\; , \;\})$. A Poisson deformation is a deformation 
of the pair of $X$ itself and the Poisson structure on it. 
When $X$ is not a compact variety, 
the usual deformation theory does not work in general 
because the tangent object $\mathbf{T}^1_X$ may possibly 
have infinite dimension, and moreover, infinitesimal or formal deformations 
do not capture actual deformations of non-compact varieties. 
On the other hand, Poisson deformations 
work very well in many important cases where $X$ is not 
a complete variety.  Denote by  ${\mathrm PD}_X$  the Poisson 
deformation functor of a symplectic variety (cf. \S 1). 
 In this paper, we shall study the 
Poisson deformation of an affine symplectic variety. 
The main result is: 
\vspace{0.2cm}

{\bf Theorem (5.1)}. {\em Let $X$ be an affine symplectic 
variety. Then the Poisson deformation functor 
$\mathrm{PD}_X$ is unobstructed.} 
\vspace{0.2cm}

A Poisson deformation of $X$ is controlled by the Poisson 
cohomology $\mathrm{HP}^2(X)$ (cf. [G-K], [Na 2]). When $X$ has only 
terminal singularities, we have $\mathrm{HP}^2(X) 
\cong H^2((X_{reg})^{an}, \mathbf{C})$, where $(X_{reg})^{an}$ is the complex 
space associated with $X_{reg}$. In that case this description enables 
us to prove that $\mathrm{PD}_X$ is unobstructed 
([Na 2], Corollary 15). But, in general, there is no 
such direct, topological description of $\mathrm{HP}^2(X)$. 
Let us explain our strategy to describe $\mathrm{HP}^2(X)$. 
As remarked, $\mathrm{HP}^2(X)$ is identified with 
$\mathrm{PD}_X(\mathbf{C}[\epsilon])$ where $\mathbf{C}[\epsilon]$ 
is the ring of dual numbers over $\mathbf{C}$. 
First, note that there is an open locus $U$ of $X$ where 
$X$ is smooth, or is locally a trivial deformation of a 
(surface) rational double point at each $p \in U$. 
Let $\Sigma$ be the singular locus of $U$. 
Note that $X \setminus U$ has codimension $\geq 4$ in 
$X$ (cf. [Ka 1]). Moreover, we have $\mathrm{PD}_X(\mathbf{C}[\epsilon]) 
\cong \mathrm{PD}_U(\mathbf{C}[\epsilon]).$ 
Put $T^1_{U^{an}} := \underline{\mathrm{Ext}}^1
(\Omega^1_{U^{an}}, \mathcal{O}_{U^{an}})$. As is well-known, 
a (local) section of $T^1_{U^{an}}$ corresponds to a 1-st 
order deformation of $U^{an}$. In \S 1, we shall construct a 
locally constant sheaf $\mathcal{H}$ of $\mathbf{C}$-modules as 
a subsheaf of $T^1_{U^{an}}$. The sheaf $\mathcal{H}$ is intrinsically 
characterized as the sheaf of germs of sections of 
$T^1_{U^{an}}$ {\em which come from Poisson deformations of} 
$U^{an}$ (cf. Lemma (1.5)). Now we have an exact sequence 
(cf. Proposition (1.11)): 
$$ 0 \to H^2(U^{an}, \mathbf{C}) \to \mathrm{PD}_U(\mathbf{C}[\epsilon])  
\to H^0(\Sigma, \mathcal{H}). $$ 
Here the first term $H^2(U^{an}, \mathbf{C})$ is  
the space of locally trivial\footnote{More exactly, this means 
that the Poisson deformations are locally trivial as usual flat 
deformations of $U^{an}$} Poisson deformations of $U$. 
By the definition of  $U$, there exists a minimal resolution 
$\pi: \tilde{U} \to U$. Let $m$ be the number of irreducible 
components of the exceptional divisor of $\pi$. Section 3 is a 
preliminary section for section 4. However, Proposition (3.2) is the  
core of the argument in \S 4.  
The main result of \S 4 is: 
\vspace{0.2cm}

{\bf Proposition (4.2)}. {\em The following equality holds:} 
$$ \dim H^0(\Sigma, \mathcal{H}) = m.$$ 

In order to prove Proposition (4.2), we need to know the monodromy 
action of $\pi_1(\Sigma)$ on $\mathcal{H}$. The idea is 
to compare two sheaves $R^2\pi^{an}_*\mathbf{C}$ and $\mathcal{H}$. 
Note that, for each point $p \in \Sigma$, the germ $(U,p)$ is isomorphic 
to the product of an ADE surface singularity $S$ and $(\mathbf{C}^{2n-2},0)$. 
Let $\tilde{S}$ be the minimal resolution of $S$.  
Then, $(R^2\pi^{an}_*\mathbf{C})_p$ is isomorphic to $H^2(\tilde{S}, \mathbf{C})$. 
A monodromy of $R^2\pi^{an}_*\mathbf{C}$ comes from a graph automorphism 
of the Dynkin diagram determined by the exceptional (-2)-curves on ${\tilde S}$. 
As is well known, $S$ is described in terms of a simple Lie algebra $\g$, and 
$H^2(\tilde{S}, \mathbf{C})$ is identified with the Cartan subalgebra $\h$ of 
$\g$; therefore, one may regard $R^2\pi^{an}_*\mathbf{C}$ as a local system 
of the $\mathbf{C}$-module $\h$ (on $\Sigma$), whose monodromy action coincides 
with the natural action of a graph automorphism on $\h$. 
On the other hand, $\mathcal{H}$ is a local system of $\h/W$, where $\h/W$ is the 
linear space obtained as the quotient of $\h$ by the Weyl group $W$ of $\g$.  
The action of a graph automorphism on $\h$ descends to an action on $\h/W$, 
which gives a monodromy action for $\mathcal{H}$. This description of the monodromy 
enables us to compute  $\dim H^0(\Sigma, \mathcal{H})$.   
 
Proposition (4.2) together with the exact sequence above 
gives an upper-bound of $\dim \mathrm{PD}_U(\mathbf{C}[\epsilon])$ in terms 
of some topological data of $X$ (or $U$). 
In \S 5, we shall prove Theorem (5.1) by using this 
upper-bound. The rough idea is the following. 
There is a natural map of functors $\mathrm{PD}_{\tilde{U}} 
\to \mathrm {PD}_U$ induced by the resolution map $\tilde{U} \to U$. 
The tangent space $\mathrm{PD}_{\tilde{U}}(\mathbf{C}[\epsilon])$ to 
$\mathrm{PD}_{\tilde{U}}$ is identified with $H^2(\tilde{U}^{an}, \mathbf{C})$. 
We have an exact sequence 
$$ 0 \to H^2(U^{an}, \mathbf{C}) \to H^2(\tilde{U}^{an}, \mathbf{C}) 
\to H^0(U^{an}, R^2\pi^{an}_*\mathbf{C}) \to 0,$$ and 
$\dim H^0(U^{an}, R^2\pi^{an}_*\mathbf{C}) = m$.  In particular, we have 
$\dim H^2(\tilde{U}^{an}, \mathbf{C}) = \dim H^2(U^{an}, \mathbf{C})  
+ m$. But this implies that $\dim \mathrm{PD}_{\tilde{U}}(\mathbf{C}[\epsilon]) 
\geq \dim \mathrm{PD}_U(\mathbf{C}[\epsilon])$. 
On the other hand, the map $\mathrm{PD}_{\tilde{U}} \to \mathrm{PD}_U$ 
has a finite closed fiber; or more exactly, the corresponding map 
$\mathrm{Spec}R_{\tilde U} \to \mathrm{Spec}R_U$ of prorepresentable hulls,  
has a finite closed fiber. 
Since $\mathrm{PD}_{\tilde{U}}$ is unobstructed, 
this implies that $\mathrm{PD}_U$ is unobstructed and 
$\dim \mathrm{PD}_{\tilde{U}}(\mathbf{C}[\epsilon]) = 
\dim \mathrm{PD}_U(\mathbf{C}[\epsilon])$. Finally, we obtain  
the unobstructedness of $\mathrm{PD}_X$ from 
that of $\mathrm{PD}_U$.      
  
Theorem (5.1) is only concerned with the formal deformations 
of $X$; but, if we impose the following condition (*), then 
the formal  universal Poisson deformation of $X$ has an 
algebraization.
\vspace{0.2cm}

(*):  $X$ has a $\mathbf{C}^*$-action with positive weights 
with a unique fixed point $0 \in X$. Moreover, $\omega$ 
is positively weighted for the action.   
\vspace{0.2cm}

We shall briefly explain how this condition (*) is used in the 
algebraization.  Let $R_X:= \lim R_X/(m_X)^{n+1}$ be the prorepresentable hull of 
$\mathrm{PD}_X$. Then the formal universal deformation $\{X_n\}$ of  $X$ 
defines an $m_X$-adic ring $A := \lim \Gamma (X_n, \mathcal{O}_{X_n})$ and 
let $\hat{A}$ be the completion of $A$ along the maximal ideal of $A$.  The rings 
$R_X$ and $\hat{A}$  both have natural $\mathbf{C}^*$-actions induced 
from the $\mathbf{C}^*$-action on $X$, and there is a $\mathbf{C}^*$-equivariant 
map $R_X \to \hat{A}$.  By taking the $\mathbf{C}^*$-subalgebras of $R_X$ and 
$\hat{A}$ generated by eigen-vectors, we get a map 
$$\mathbf{C}[x_1, ..., x_d] \to S$$ from a polynomial ring to a $\mathbf{C}$-algebra 
of finite type. We also have a Poisson structure on $S$ over $\mathbf{C}[x_1, ..., x_d]$ 
by the second condition of (*).  As a consequence,  
there is an affine space $\mathbf{A}^d$ whose 
completion at the origin coincides with $\mathrm{Spec}(R_X)$ in such 
a way that the formal universal Poisson deformation over $\mathrm{Spec}(R_X)$ 
is algebraized to a $\mathbf{C}^*$-equivariant map 
$$ \mathcal{X} \to \mathbf{A}^d. $$   
Now, by using the minimal model theory due to Birkar-Cascini-Hacon-McKernan 
[BCHM], one can study the general fiber of $\mathcal{X} \to \mathbf{A}^d$.  
According to [BCHM], we can take a crepant partial resolution 
$\pi: Y \to X$ in such a way that $Y$ has only $\mathbf{Q}$-factorial 
terminal singularities. This $Y$ is called a {\bf Q}-{\em factorial terminalization} of $X$.  
In our case, $Y$ is a symplectic variety and 
the $\mathbf{C}^*$-action on $X$ uniquely extends to that on $Y$. 
Since $Y$ has only terminal singularities, it is relatively easy to 
show that the Poisson deformation functor $\mathrm{PD}_Y$ is 
unobstructed. Moreover, the formal universal Poisson deformation of 
$Y$ has an algebraization over an affine space $\mathbf{A}^d$: 
$$  \mathcal{Y} \to \mathbf{A}^d. $$ 
There is a $\mathbf{C}^*$-equivariant commutative diagram 

\begin{equation} 
\begin{CD} 
\mathcal{Y} @>>> \mathcal{X}  \\ 
@VVV @VVV \\ 
\mathbf{A}^d @>{\psi}>> \mathbf{A}^d    
\end{CD} 
\end{equation} 

By Theorem (5.5), (a): $\psi$ is a finite surjective map, 
(b):  $\mathcal{Y} \to \mathbf{A}^d$ is a 
locally trivial deformation of $Y$, and 
(c): the induced map 
$\mathcal{Y}_t \to \mathcal{X}_{\psi(t)}$ is an isomorphism 
for a general point $t \in \mathbf{A}^d$. 
As an application of Theorem (5.5), we have       
\vspace{0.2cm}

{\bf Corollary (5.6)}: {\em Let $(X, \omega)$ be an affine symplectic 
variety with the property (*).  Then the following are equivalent.} 
\vspace{0.15cm} 

(1) {\em $X$ has a crepant projective resolution.} 
\vspace{0.12cm}

(2) {\em $X$ has a smoothing by a Poisson deformation.} 
\vspace{0.2cm}

{\bf Example} (i) Let $O \subset \g$ be a nilpotent orbit of a complex 
simple Lie algebra. Let $\tilde{O}$ be the normalization of the 
closure $\bar{O}$ of $O$ in $\g$. Then $\tilde{O}$ is an affine symplectic 
variety with the Kostant-Kirillov 2-form $\omega$ on $O$.   
Let $G$ be a complex algebraic group with $Lie(G) = \g$. 
By [Fu], $\tilde{O}$ has a crepant projective resolution if and only if 
$O$ is a Richardson orbit (cf. [C-M]) and there is a parabolic subgroup 
$P$ of $G$ such that its Springer map $T^*(G/P) \to \tilde{O}$ is birational. 
In this case, every crepant resolution of $\tilde{O}$ is actually obtained as 
a Springer map for some $P$. If $\tilde{O}$ has a crepant resolution, 
$\tilde{O}$ has a smoothing by a Poisson deformation. The smoothing of 
$\tilde{O}$ is isomorphic to the affine variety $G/L$, where $L$ is the 
Levi subgroup of $P$. Conversely, if $\tilde{O}$ has a smoothing by a Poisson deformation, 
then the smoothing always has this form.   

(ii) In general, $\tilde{O}$ has no crepant resolutions. But   
a suitable generalized Springer map gives a {\bf Q}-factorial terminalization 
of $\tilde{O}$ by [Na 4] and [Fu 2].  
More explicitly, there 
is a parabolic subalgebra $\p$ with Levi decomposition $\p = \n \oplus \l$ and 
a nilpotent orbit $O'$ in $\l$ so that the generalized Springer map 
$G \times^P (\n + \bar{O'}) \to \tilde{O}$ is a crepant, birational map, and  
the normalization of  $G \times^P (\n + \bar{O'})$ is a {\bf Q}-factorial terminalization 
of $\tilde{O}$.  By a Poisson deformation, $\tilde{O}$ deforms to the normalization of 
$G \times^L \bar{O'}$. Here $G \times^L \bar{O'}$ is a fiber bundle over $G/L$ with 
a typical fiber $\bar{O'}$, and its normalization can be written as $G \times^L \tilde{O'}$ 
with the normalization $\tilde{O'}$ of $\bar{O'}$.  
\vspace{0.15cm}

\section{Local system associated with a symplectic 
variety} 
 
(1.1) A {\em symplectic variety} $(X, \omega)$ is a pair of a normal algebraic variety $X$ 
defined over {\bf C}  
and a symplectic 2-form $\omega$ on the regular part $X_{reg}$ of $X$ such 
that, for any resolution $\mu: \tilde{X} \to X$, the 2-form $\omega$ on $\mu^{-1}(X_{reg})$ 
extends to a closed regular 2-form on $\tilde{X}$.  We also have a similar notion of a symplectic variety in 
the complex analytic category (eg. the germ of a normal complex space, a holomorphically 
convex, normal, complex space).  For an algebraic variety $X$ over {\bf C}, we denote by 
$X^{an}$ the associated complex space. Note that if $(X, \omega)$ is a symplectic variety, then 
$X^{an}$ is naturally a symplectic variety in the complex analytic category.  
A symplectic variety $X$ (resp. $X^{an}$) has rational Gorenstein singularities.  
The symplectic 2-form $\omega$ defines a bivector $\Theta \in \wedge^2 \Theta_{X_{reg}}$ 
by the identification $\Omega^2_{X_{reg}} \cong \wedge^2 \Theta_{X_{reg}}$ by $\omega$. 
Define a Poisson structure  
$\{ \;, \; \}$ on $X_{reg}$ by $\{f,g\} := \Theta(df \wedge dg)$. Since $X$ is normal, the  
Poisson structure on $X_{reg}$ uniquely extends to a Poisson structure on $X$. 
Here, we recall 
the definition of a Poisson scheme or a Poisson complex space.   
 
\vspace{0.2cm}
{\bf Definition}. Let $T$ be a scheme (resp. complex space).  
Let $X$ be a scheme  (resp.  
complex space) over $T$. 
Then $(X, \{\; , \;\})$ is a Poisson scheme (resp. a Poisson space) over $T$ if 
$\{\; , \;\}$ is an $\mathcal{O}_T$-linear map: 
$$ \{\; , \;\}: \wedge^2_{\mathcal{O}_T}\mathcal{O}_X \to 
\mathcal{O}_X$$ 
such that, for $a,b,c \in \mathcal{O}_X$, 
\begin{enumerate}
\item $\{a, \{b,c\}\} + \{b,\{c,a\}\} + \{c,\{a,b\}\} = 0$ 
\item $\{a,bc\} = \{a,b\}c + \{a,c\}b.$ 
\end{enumerate}

Let $(X, \{\;, \;\})$ be a Poisson scheme (resp. Poisson space) over $\mathbf{C}$. 
Let $S$ be a local Artinian $\mathbf{C}$-algebra with 
$S/m_S = \mathbf{C}$.  Let $T$ be the affine scheme (resp. complex space) whose 
coordinate ring is $S$.  A Poisson deformation of $(X, \{\;, \;\})$ 
over $S$ is a Poisson scheme (resp. Poisson complex space) over $T$: $(\mathcal{X}, \{\;, \;\}_T)$ 
such that $\mathcal{X}$ is flat over $T$, 
$\mathcal{X} \times_T \mathrm{Spec}(\mathbf{C}) 
\cong X$, and the Poisson structure $\{\; , \;\}_T$ induces 
the original Poisson structure $\{\; , \;\}$ over the closed 
fiber $X$.   
We define 
$\mathrm{PD}_X(S)$ to be the set of equivalence classes 
of the pairs of Poisson deformations $\mathcal{X}$ of $X$ over 
$\mathrm{Spec}(S)$ and Poisson isomorphisms 
$\phi: \mathcal{X}\times_{\mathrm{Spec}(S)}\mathrm{Spec}(\mathbf{C}) 
\cong X$. Here $(\mathcal{X}, \phi)$ and 
$(\mathcal{X}', \phi')$ are equivalent if there is a Poisson 
isomorphism $\varphi: \mathcal{X} \cong \mathcal{X}'$ over 
$\mathrm{Spec}(S)$ which induces the identity map of $X$ over 
$\mathrm{Spec}(\mathbf{C})$ via $\phi$ and $\phi'$. 
We define the {\em Poisson deformation functor}:  
$$\mathrm{PD}_{(X, \{\;, \;\})}: (\mathrm{Art})_{\mathbf{C}} 
\to (\mathrm{Set})$$ from the category of local Artin $\mathbf{C}$-algebras 
with residue field $\mathbf{C}$ to the category of sets.  
Let $\mathbf{C}[\epsilon]$ be the ring of dual numbers over $\mathbf{C}$.  
Then the set $\mathrm{PD}_X(\mathbf{C}[\epsilon])$ has a structure of the {\bf C}-vector space, and 
it is called the tangent space of $\mathrm{PD}_X$. 
A Poisson deformation of $X$ over $\mathrm{Spec}\mathbf{C}[\epsilon]$ is 
particularly called a {\em 1-st order} Poisson deformation of $X$. 
It is easy to see that $\mathrm{PD}_{(X, \{\;, \;\})}$ satisfies 
the Schlessinger's conditions ([Sch]) except that possibly $\dim \mathrm{PD}_{(X, \{\;, \;\})}
(\mathbf{C}[\epsilon]) = \infty $.
For details on Poisson deformations, see [G-K], [Na 2].  

(1.2) Let $(S, 0)$ be the germ of a rational double point of 
dimension $2$.  More explicitly, 
$$ S :=\{(x,y,z) \in \mathbf{C}^3; f(x,y,z) = 0\}, $$ 
where 
$$ f(x,y,z) = xy + z^{r+1}, $$ 
$$ f(x,y,z) = x^2 +   y^2z + z^{r-1}, $$  
$$ f(x,y,z) = x^2 + y^3 + z^4, $$ 
$$ f(x,y,z) = x^2 + y^3 + yz^3, $$ or 
$$ f(x,y,z) = x^2 + y^3 + z^5$$   
according as $S$ is of type $A_r$, $D_r$ ($r \geq 4$) 
$E_6$, $E_7$ or $E_8$. 
We put 
$$ \omega_S := res(dx \wedge dy \wedge dz/f). $$ 
Then $\omega_S$ is a symplectic 2-form on 
$S - \{0\}$ and $(S,0)$ becomes a symplectic variety.   
Let us denote by $\omega_{\mathbf{C}^{2m}}$ the 
canonical symplectic form on $\mathbf{C}^{2m}$ : 
$$ ds_1 \wedge dt_1 +  ... + ds_m 
\wedge dt_m. $$ 

Let $(X, \omega)$ be a symplectic variety of dimension $2n$ whose  
singularities are (analytically) locally isomorphic to 
$(S,0) \times (\mathbf{C}^{2n-2},0)$.  Let $\Sigma$ be the 
singular locus of $X$. \vspace{0.2cm} 

{\bf Lemma (1.3)}  {\em For any $p \in \Sigma$, there are  
an open neighborhood $U \subset X^{an}$ of $p$ and an  
open immersion  
$$ \phi : U \to S \times \mathbf{C}^{2n-2} $$ 
such that $\omega\vert_U = \phi^*((p_1)^*\omega_S +  
(p_2)^*\omega_{\mathbf{C}^{2n-2}})$, where $p_i$ are 
$i$-th projections of $S \times \mathbf{C}^{2n-2}$.} 
\vspace{0.2cm} 

{\em Proof}.  
Let $\omega_1$ be an arbitrary symplectic 2-form on 
the regular locus of $(S,0) \times (\mathbf{C}^{2n-2},0)$. 
On the other hand, we put $$\omega_0 := 
(p_1)^*\omega_S +  (p_2)^*\omega_{\mathbf{C}^{2n-2}}.$$  
The singularity $(S,0)$ can be written as $(\mathbf{C}^2,0)/G$ 
with a finite subgroup $G \subset SL(2, \mathbf{C})$. Let 
$\pi: (\mathbf{C}^2, 0) \to (S,0)$ be the quotient map. The finite 
group $G$ acts on  $(\mathbf{C}^2,0) \times (\mathbf{C}^{2n-2},0)$ in such a way 
that it acts on the second factor trivially. Then one has  
the quotient map $$ \pi \times id: (\mathbf{C}^2,0) \times (\mathbf{C}^{2n-2},0) 
\to (S,0) \times (\mathbf{C}^{2n-2},0).$$ We put 
$$ \tilde{\omega}_i := (\pi \times id)^*{\omega_i}$$ for $i = 0,1$.  
Then $\tilde{\omega}_i$ are $G$-invariant symplectic 2-forms on 
$(\mathbf{C}^2,0) \times (\mathbf{C}^{2n-2},0)$. 
We shall prove that there is a $G$-equivarinat automorphism 
$\tilde{\varphi}$ of $(\mathbf{C}^2,0) \times (\mathbf{C}^{2n-2},0)$ 
such that $\tilde{\varphi}^*\tilde{\omega}_1 = \tilde{\omega}_0$.  The basic idea 
of the following arguments is due to [Mo]. 
Let $(x,y)$ be the coordinates of $(\mathbf{C}^2,0)$ and 
let $(s_1, ..., s_{n-1}, t_1, ..., t_{n-1})$ be the coordinates of 
$(\mathbf{C}^{2n-2},0)$.  
The symplectic 2-forms $\tilde{\omega}_0$ and $\tilde{\omega}_1$ restrict respectively  
to give 2-forms $\tilde{\omega}_0(\mathbf{0})$ and $\tilde{\omega}_1(\mathbf{0})$ on 
the tangent space $T_{\mathbf{C}^{2n}, \mathbf{0}}$ at the origin $\mathbf{0} \in \mathbf{C}^{2n}$. 
By the definition of $\tilde{\omega}_0$,   
$$\tilde{\omega}_0(\mathbf{0}) = adx \wedge dy + \Sigma ds_i \wedge dt_i$$ with 
some $a \in \mathbf{C}^*$.  
Next write $\tilde{\omega}_1(\mathbf{0})$ by using $dx$, $dy$, $ds_i$ and $dt_j$. 
We may assume that $G$ contains a diagonal matrix   
$$\left( \begin{array}{cc}
\zeta & 0 \\ 
0 & \zeta^{-1} 
\end{array}\right)$$ where $\zeta$ is a primitive $l$-th root of unity with 
some $l > 1$. Since $\tilde{\omega}_1$ is $G$-invariant,  
$\tilde{\omega}_1(\mathbf{0})$ does not contain the terms $dx \wedge ds_i$, 
$dx \wedge dt_j$, $dy \wedge ds_i$ or $dy \wedge dt_j$.  
One can choose a scalar multiplication $c: ({\mathbf C}^2,0) \to ({\mathbf C}^2,0)$ 
($(x,y) \to (cx,cy)$) and a linear automorphism $\sigma : ({\mathbf C}^{2n-2},0) 
\to ({\mathbf C}^{2n-2},0)$ so that $\tilde{\omega}_2 := 
(c \times \sigma)^*(\tilde{\omega}_1)$ satisfies 
$$ \tilde{\omega}_2(\mathbf{0}) = adx \wedge dy  + \Sigma ds_i \wedge dt_i.$$     
Note that $$ \tilde{\omega}_0(\mathbf{0}) = \tilde{\omega}_2(\mathbf{0}).$$ 
Since $c \times \sigma$ is $G$-equivariant, $\tilde{\omega}_2$ is a $G$-invariant 
symplectic 2-form.  
For $\tau \in {\mathbf R}$, define 
$$ \omega({\tau}) := (1-\tau) \tilde{\omega}_0 + \tau \tilde{\omega}_2.$$ 
We put $$u := d\omega({\tau})/d\tau.$$ 
Since $S \times {\mathbf C}^{2n-2}$ has only quotient singularities, the complex 
$((\pi \times id)^G_*\Omega^{\cdot}_{{\mathbf C}^2 \times {\mathbf  C}^{2n-2}}, d)$ is a 
resolution of the constant sheaf ${\mathbf C}$ on
$S \times {\mathbf C}^{2n-2}$. Note that $u$ is a section 
of $(\pi \times id)^G_*\Omega^2_{{\mathbf C}^2 \times {\mathbf  C}^{2n-2}}$. Moreover, 
$u$ is d-closed. Therefore, one can write $u = dv$ with a $G$-invariant 1-form $v$. 
Moreover $v$ can be chosen such that $v(\mathbf{0}) = 0$.  
Define a vector field $X_{\tau}$ on $(\mathbf{C}^{2n}, 0)$ by 
$$i_{X_{\tau}}\omega(\tau) = -v.$$ 
Since $\omega(\tau)$ is $d$-closed, we have $$L_{X_{\tau}}\omega(\tau) = -u$$ 
where $L_{X_{\tau}}\omega(\tau)$ is the 
Lie derivative of $\omega(\tau)$ along $X_{\tau}$.   
If we take a sufficiently small open subset $V$ of $\mathbf{0} \in \mathbf{C}^{2n}$, then the vector fields  
$\{X_{\tau}\}_{0 \le \tau \le 1}$ 
define a family of open immersions   
$\varphi_{\tau} : V \to {\mathbf C}^{2n}$ via 
$$ d\varphi_{\tau}/d\tau  = X_{\tau}(\varphi_{\tau}), \;\; \varphi_0 = id.$$ 
Since all $\varphi_{\tau}$ fix the origin and $X_{\tau}$ are all $G$-invariant, 
$\varphi_{\tau}$ induce $G$-equivariant automorphisms of $({\mathbf C}^{2n},0)$. 
By the definition of $X_{\tau}$, we have 
$(\varphi_{\tau})^*\omega(\tau) = \omega(0)$. In particular, 
$(\varphi_1)^*\tilde{\omega}_2 = \tilde{\omega}_0$.  We put $$\tilde{\varphi} := (\varphi_1)\circ 
(c \times \sigma).$$  
The $G$-equivariant automorphism $\tilde{\varphi}$ of $(\mathbf{C}^{2n},0)$ descends 
to an automorphism $\varphi$ of $(S,0) \times (\mathbf{C}^{2n-2},0)$ so that 
$\varphi^*\omega_1 = \omega_0$.  
Q.E.D. 
\vspace{0.2cm}

We cover the singular locus $\Sigma$ by a family of open sets $\{U_{\alpha}\}$ 
of $X^{an}$ in such a way that each $U_{\alpha}$ 
admits an open immersion  $\phi_{\alpha}$ as in Lemma (1.3). 
In the remainder, we call such a 
covering $\{U_{\alpha}\}$ {\em admissible}.

(1.4)  Let $(X, \omega)$ be the same as above. 
Denote by $T^1_{X^{an}}$ the analytic coherent sheaf  
$\underline{\mathrm{Ext}}^1(\Omega^1_{X^{an}}, 
\mathcal{O}_{X^{an}})$. Note that the sheaf $T^1_{X^{an}}$ is the sheafication 
of the presheaf associating to each open set $V \subset X^{an}$ the 
$\mathbf{C}$-vector space of the isomorphic classes of 1-st order deformations of $V$. 
Let us consider the presheaf 
on $X^{an}$ which associates to each open set $V$ the $\mathbf{C}$-vector space of the isomorphic 
classes of 1-st order {\em Poisson} deformation. Denote by $PT^1_{X^{an}}$ the sheafication 
of this presheaf. Note that both sheaves $T^1_{X^{an}}$ and $PT^1_{X^{an}}$ have  
support on $\Sigma$. 
One has a natural map $$PT^1_{X^{an}} \to T^1_{X^{an}}$$ of sheaves of 
$\mathbf{C}$-modules by forgetting the Poisson structure. 
Define a subsheaf $\mathcal{H}$ of $T^1_{X^{an}}$ as the image of this map.   
\vspace{0.2cm}
  
{\bf Lemma (1.5)} {\em $\mathcal{H}$ is a locally constant $\mathbf{C}$-module 
over  $\Sigma$. }   
\vspace{0.2cm}

{\em Proof}. 
Take an admissible covering $\{U_{\alpha}\}$.  For each $\alpha$, 
$$T^1_{U_{\alpha}} = (p_1\circ \phi_{\alpha} )^*T^1_S.$$   
We put $$ H_{\alpha} := (p_1\circ \phi_{\alpha})^{-1}T^1_S.$$ 
Note that $H_{\alpha}$ is a constant $\mathbf{C}$-module on $U_{\alpha} \cap \Sigma$. 
We shall prove that $\mathcal{H}\vert_{U_{\alpha}} = H_{\alpha}$.    
In fact, let $\mathcal{U}_{\alpha} \to \mathrm{Spec}\mathbf{C}[\epsilon]$ 
be a 1-st order Poisson deformation of $U_{\alpha}$. 
Let $0 \in U_{\alpha}$ be the point which corresponds to 
$(0,0) \in S \times \mathbf{C}^{2n-2}$ via $\phi_{\alpha}$. 
By applying the second statement of the next Lemma (1.6) to 
$\hat{\mathcal{O}}_{U_{\alpha},0}$ and 
$\hat{\mathcal{O}}_{\mathcal{U}_{\alpha}, 0}$, we conclude 
that $(\mathcal{U}_{\alpha},0) \cong 
(\mathcal{S}, 0) \times (\mathbf{C}^{2n-2}, 0)$, 
where $\mathcal{S}$ is a 1-st order deformation of $S$.   
Conversely, a 1-st order deformation of this form 
always comes from a Poisson deformation of $U_{\alpha}$. Q.E.D. 
\vspace{0.2cm}

{\bf Lemma (1.6)}. {\em Let $S := \{f(x,y,z) = 0\} \subset \mathbf{C}^3$ be an 
isolated hypersurface singularity which admits a Poisson structure, and 
let $(\mathbf{C}^{2n-2},0)$ be a symplectic manifold with the standard 
symplectic structure. Put $V := (S,0) \times (\mathbf{C}^{2n-2},0)$ and 
introduce the product Poisson structure on $V$. Assume that 
$\mathcal{V} \to \mathrm{Spec}\;\mathbf{C}[\epsilon]$ 
is a 1-st order Poisson deformation of $V$. Then $$\mathcal{V} \cong  
(\mathcal{S},0) \times (\mathbf{C}^{2n-2},0)$$ as a flat deformation. Here 
$(\mathcal{S},0)$ is a 1-st order flat deformation of $(S,0)$.}  

{\em Proof}. We denote by $\mathbf{s} = (s_1, ..., s_{2n-2})$ the coordinates 
of $\mathbf{C}^{2n-2}$. 
Let $f_1, ..., f_{\tau} \in \mathbf{C}\{x,y,z\}$ be 
the representatives of a basis of $\mathbf{C}\{x,y,z\}/(f, f_x, f_y, f_z)$. 
The 1-st order deformation $\mathcal{V}$ can be written as 
$$ f(x,y,z) + \epsilon (f_1(x,y,z)g_1(\mathbf{s}) + ... 
+ f_{\tau}(x,y,z)g_{\tau}(\mathbf{s})) = 0.$$ 
We prove that $g_i$ are all constants. Let $\{\;, \;\}$ be the Poisson structure 
on $V$. By the definition, we have $$\{x, s_i\} = \{y,s_i\} = \{z, s_i\} = 0$$ in 
$\mathcal{O}_{V,0}$. Let $\{\;, \;\}'$ be the Poisson structure on $\mathcal{V}$ 
extending the Poisson structure $\{\;, \;\}$. Then we have 
$$ \{x, s_i\}' = \epsilon \alpha_i, \; \{y, s_i\}' = \epsilon \beta_i, \; 
\{z, s_i\}' = \epsilon \gamma_i,$$ for some elements $\alpha_i$, $\beta_i$ and 
$\gamma_i$ in $\mathcal{O}_{V,0}$. 
Since $f + \epsilon (f_1g_1 + ... + f_{\tau}g_{\tau}) = 0$ in $\mathcal{O}_{\mathcal{V},0}$, 
we must have 
$$ \{f + \epsilon (f_1g_1 + ... + f_{\tau}g_{\tau}), s_i\}' = 0$$ in 
$\mathcal{O}_{\mathcal{V},0}$. 
By calculating the left-hand side, one has  
$$ f_x\{x,s_i\}' + f_y\{y, s_i\}' + f_z\{z, s_i\}' + 
\epsilon( \Sigma_{1 \le j \le \tau} f_j\{g_j, s_i\} + 
\Sigma_{1 \le j \le \tau} g_j\{f_j, s_i\}) = 0$$ 
Recall that 
$$\{s_1, s_2\} = \{s_3, s_4\} = ... = \{s_{2n-3}, s_{2n-2}\} = 1,$$ 
and $\{s_k, s_l\} = 0$ for other $k < l$. Moreover, 
note that $\{f_j, s_i\} = 0$. Assume that $i$ is odd, then one has 
$$ \epsilon (f_x \alpha_i + f_y\beta_i + f_z\gamma_i 
+ \Sigma_{1 \le i \le \tau}f_j\cdot (g_j)_{s_{i+1}}) = 0.$$    
This implies that 
$$ f_x \alpha_i + f_y\beta_i + f_z\gamma_i 
+ \Sigma_{1 \le i \le \tau}f_j\cdot (g_j)_{s_{i+1}} = 0$$ in 
$\mathcal{O}_{V,0}$. Note that $\mathcal{O}_{V,0} 
= \mathbf{C}\{x,y,z, \mathbf{s}\}/(f)$. Let us consider the 
equation in $\mathbf{C}\{x,y,z, \mathbf{s}\}/(f, f_x, f_y, f_z)$. 
Then we have $$\Sigma_{1 \le j \le \tau} f_j\cdot(g_j)_{s_{i+1}} = 0.$$
This implies that $(g_j)_{s_{i+1}} = 0$ in $\mathbf{C}\{\mathbf{s}\}$ 
for all $j$. When $i$ is even, a similar argument shows that 
$(g_j)_{s_{i-1}} = 0$ for all $j$. As a consequence, $g_j$ are 
constants for all $j$.  Q.E.D.    
\vspace{0.2cm} 

(1.7) {\bf Monodromy of $\mathcal{H}$} 

Let $\gamma$ be a closed loop in $\Sigma$ 
starting from $p \in \Sigma$. 
We shall describe the monodromy of $\mathcal{H}$
along $\gamma$ in terms of a certain symplectic automorphism 
of the germ $(X^{an},p)$.  
In order to do this, we take a sequence of admissible 
open sets of $X^{an}$: $U_1$, ..., $U_k$, $U_{k+1} := U_1$ 
in such a way that $p \in U_1$, $\gamma \subset 
\cup U_i$, $U_i \cap U_{i+1} \cap \gamma \ne \emptyset$ for 
$i = 1, ..., k$.  
Put $p_1 := p$ and 
choose a point $p_i \in U_i \cap U_{i+1} \cap \gamma$ for each 
$i \geq 2$.  Let $\phi_i : U_i \to S \times \mathbf{C}^{2n-2}$ be 
the symplectic open immersion associated with the admissible open 
subset $U_i$.  Since $\mathcal{H}$ is a locally constant 
$\mathbf{C}$-module by (1.5),  an element of $\mathcal{H}_{p_i}$ uniquely extends to a 
section of $\mathcal{H}$ over $U_i$. Since $p_{i-1} \in U_i$, this section 
restricts to give an element of $\mathcal{H}_{p_{i-1}}$. In this way, we have 
an identification $$m_i : \mathcal{H}_{p_{i-1}} \cong \mathcal{H}_{p_i}$$ for each $i$. 
The monodromy transformation $m_{\gamma}$ is the composite of $m_i$'s: 
$$ m_{\gamma} = m_{k+1} \circ ... \circ m_{2}. $$ One can describe each $m_i$ in terms 
of certain symplectic isomorphisms as explained below.  
Since $U_i$ contains $p_i$,  the germ $(X^{an}, p_i)$ is identified with 
$(S \times \mathbf{C}^{2n-2}, \phi_i(p_i))$ 
by $\phi_i$. On the other hand, since $U_i$ contains $p_{i-1}$, the germ 
$(X^{an}, p_{i-1})$ is identified with $(S \times \mathbf{C}^{2n-2}, \phi_i(p_{i-1}))$. 
Note that $\phi_i(p_i) = (0, *) \in S \times \mathbf{C}^{2n-2}$ and 
$\phi_i(p_{i-1}) = (0, **) \in S \times \mathbf{C}^{2n-2}$ for some points 
$*, ** \in \mathbf{C}^{2n-2}$ because $p_i, p_{i-1} \in \gamma$.  
Denote by $\sigma_i : \mathbf{C}^{2n-2} \to \mathbf{C}^{2n-2}$ the 
translation map such that $\sigma_i(*) = **$.  Then, by the automorphism  
$id \times \sigma_i $ of  $S \times \mathbf{C}^{2n-2}$, two germs 
$(S \times \mathbf{C}^{2n-2}, \phi_i(p_i))$ and 
$(S \times \mathbf{C}^{2n-2}, \phi_i(p_{i-1}))$ are identified. 
As a consequence, two germs $(X^{an}, p_{i-1})$ and $(X^{an}, p_i)$ have 
been identified.  By definition, this identification preserves the 
natural symplectic forms on $(X^{an}, p_{i-1})$ and $(X^{an}, p_i)$. 
The symplectic isomorphism $(X^{an}, p_{i-1}) \cong (X^{an}, p_i)$ 
determines an isomorphism $\mathcal{H}_{p_{i-1}} \cong \mathcal{H}_{p_i}$. 
It is easy to see that this isomorphism coincides with $m_i$ defined above. 
Note that the symplectic automorphism depends on the choice of 
$\phi_i$, but $m_i$ is independent of it. 
Now the sequence of identifications $(X^{an}, p_1) \cong (X^{an}, p_2)$, 
$(X^{an}, p_2) \cong (X^{an}, p_3)$, ..., $(X^{an}, p_k) \cong (X^{an}, p_1)$ 
finally defines a symplectic automorphism 
$$ {i}_{\gamma} : (X^{an}, p) \cong (X^{an}, p). $$ 
The map ${i}_{\gamma}$ induces an automorphism of $\mathcal{H}_p$, 
which coincides with $m_{\gamma}$ because $m_{\gamma} = m_{k+1} \circ ... \circ m_{2}$. 
Although $i_{\gamma}$ depends on the choices of $\phi_i$'s, $m_{\gamma}$ is 
independent of them by the definition. 
\vspace{0.2cm}

(1.8) In the above, we only considered a symplectic variety 
whose singularities are locally isomorphic to 
$(S,0) \times ({\mathbf C}^{2n-2},0)$. From now on, we will treat  
a general symplectic variety $(X, \omega)$.   
Let $U \subset X$ be the locus where 
$X$ is smooth, or is locally a trivial deformation of a 
(surface) rational double point.  
Put $\Sigma := \mathrm{Sing}(U)$. 
As an open set of $X$, $U$ naturally becomes a 
Poisson scheme. Since $X \setminus U$ has codimension 
at least $4$ in $X$ ([Ka 1]), one can prove in the same 
way as [Na 2, Proposition 13] that 
$$ \mathrm{PD}_X(\mathbf{C}[\epsilon]) \cong 
\mathrm{PD}_U(\mathbf{C}[\epsilon]). $$ 
Let $\mathrm{PD}_{lt, U}$ be the locally 
trivial Poisson deformation functor of $U$. 
More exactly, $\mathrm{PD}_{lt,U}$ is the 
subfunctor of $PD_U$ corresponding to the Poisson 
deformations of $U$ which are locally trivial 
as flat deformations of $U^{an}$ (after forgetting 
Poisson structure).  We shall insert a lemma here, which will be used in the 
proof of Proposition (1.11). 
\vspace{0.2cm}

{\bf Lemma (1.9)} {\em Let $X$ be an affine symplectic variety let 
$j: X_{reg} \to X$ be the open immersion of the regular part $X_{reg}$ 
into $X$. Then $$\mathrm{PD}_{lt,X}(\mathbf{C}[\epsilon]) = 
\mathbf{H}^2(\Gamma(X, j_*(\wedge^{\geq 1}\Theta_{X^{reg}})),$$ where    
$(\wedge^{\geq 1}\Theta_{X^{reg}}, \delta)$ is the 
Lichnerowicz-Poisson complex for $X_{reg}$ (cf. [Na 2, \S 2]). } 
\vspace{0.2cm}

{\em Proof}. The 2-nd 
cohomology $\mathbf{H}^2(\Gamma(X_{reg}, \wedge^{\geq 1}\Theta_{X_{reg}}))$ 
describes the equivalence classes of the extension of the Poisson structure $\{\;, \; \}$ on $X_{reg}$ 
to that on $X_{reg} \times \mathrm{Spec}\mathbf{C}[\epsilon] \to \mathrm{Spec}\mathbf{C}[\epsilon]$. 
In fact, for $\psi \in \Gamma (X_{reg}, \wedge^2 \Theta_{X_{reg}})$, we define a Poisson structure 
$\{\:, \:\}_{\epsilon}$ on $\mathcal{O}_{X_{reg}} \oplus \epsilon \mathcal{O}_{X_{reg}}$ 
by $$ \{f + \epsilon f', g + \epsilon g'\}_{\epsilon} := \{f, g\} + \epsilon (\psi (df \wedge dg) 
+ \{f, g'\} + \{f', g\}).$$ Then this bracket is a Poisson bracket if and only if 
$\delta (\psi) = 0$. On the other hand, an element $\theta \in \Gamma (X_{reg}, \Theta_{X_{reg}})$ 
corresponds to an automorphism $\varphi_{\theta}$ of $X_{reg} \times \mathrm{Spec}\mathbf{C}[\epsilon]$ over 
$\mathrm{Spec}\mathbf{C}[\epsilon]$ which restricts to give the identity map of the closed 
fiber $X_{reg}$. Let $\{\:, \:\}_{\epsilon}$ and $\{\:, \:\}'_{\epsilon}$ be the Poisson structures 
determined respectively by $\psi \in \Gamma(X_{reg}, \wedge^2\Theta_{X_{reg}})$ and 
$\psi' \in \Gamma(X_{reg}, \wedge^2\Theta_{X_{reg}})$. Then the two Poisson structures are  
equivalent under $\varphi_{\theta}$ if and only if $\psi - \psi' = \delta (\theta)$. 
For an affine variety $X$, a locally trivial infinitesimal deformation is nothing but a trivial 
infinitesimal deformation because $H^1(X, \Theta_X) = 0$. The original Poisson structure on 
$X$ restricts to give a Poisson structure on $X_{reg}$. As seen above, its extension to          
$X_{reg} \times \mathrm{Spec}\mathbf{C}[\epsilon]$ is classified by   
$\mathbf{H}^2(\Gamma(X_{reg}, \wedge^{\geq 1}\Theta_{X_{reg}}))$. Each Poisson structure on 
$X_{reg} \times \mathrm{Spec}\mathbf{C}[\epsilon]$ can extend uniquely to that 
on $X \times \mathrm{Spec}\mathbf{C}[\epsilon]$. 
\vspace{0.2cm}

{\bf Remark (1.10)}. By the same argument as [Na 2], Proposition 8, 
one can prove that, for a (non-affine) symplectic variety $X$, 
$$\mathrm{PD}_{lt,X}(\mathbf{C}[\epsilon]) =  
\mathbf{H}^2(X, j_*(\wedge^{\geq 1}\Theta_{X_{reg}})),$$ 
where $\mathbf{H}^2$ is the 2-nd hypercohomology. 
\vspace{0.2cm}

Let us return to the original situation in (1.8).  
Let $\mathcal{H} \subset 
T^1_{U^{an}}$ be the local constant $\mathbf{C}$-modules 
over $\Sigma$. 
We have an exact sequence of 
$\mathbf{C}$-vector spaces: 
$$ 0 \to \mathrm{PD}_{lt, U}(\mathbf{C}[\epsilon]) 
\to \mathrm{PD}_U(\mathbf{C}[\epsilon]) \to 
H^0(\Sigma, \mathcal{H}). $$  The following proposition 
shows that the tangent space of the Poisson deformation functor 
of an affine symplectic variety is finite dimensional.  
\vspace{0.2cm}
  
{\bf Proposition (1.11)}. {\em Assume that 
$X$ is an affine symplectic variety. Then 
$$ \mathrm{PD}_{lt, U}(\mathbf{C}[\epsilon]) 
\cong H^2(U^{an}, \mathbf{C}).$$ In particular, 
$\dim \mathrm{PD}_X(\mathbf{C}[\epsilon]) < \infty$.} 
\vspace{0.2cm}

{\em Proof}. Let $U^0$ be the smooth part of $U$ 
and let $j: U^0 \to U$ be the inclusion map. 
Let $(\wedge^{\geq 1}\Theta_{U^0}, \delta)$ be 
the Lichnerowicz-Poisson complex for $U^0$. By Remark (1.10), one has  
$$ \mathrm{PD}_{lt,U}(\mathbf{C}[\epsilon]) 
\cong \mathbf{H}^2(U, j_*(\wedge^{\geq 1}\Theta_{U_0})).$$ 
By the symplectic form $\omega$, the complex 
$(j_*(\wedge^{\geq 1}\Theta_{U_0}), \delta )$ 
is identified with 
$\{j_*(\wedge^{\geq 1}\Omega^1_{U_0}), d)$ 
(cf. [Na 2, Proposition 9]). The latter complex is 
the truncated {\em de Rham complex for a $V$-manifold $U$} 
$(\tilde{\Omega}^{\geq 1}_U, d)$ (cf. [St]). 
Let us consider the distinguished triangle 
$$ \tilde{\Omega}^{\geq 1}_U \to 
\tilde{\Omega}^{\cdot}_U \to \mathcal{O}_U 
\to \tilde{\Omega}^{\geq 1}_U[1].$$ 
We have an exact sequence 
$$H^1(\mathcal{O}_U) \to \mathbf{H}^2(\tilde{\Omega}^{\geq 1}_U) 
\to \mathbf{H}^2(\tilde{\Omega}^{\cdot}_U) \to H^2(\mathcal{O}_U).$$
Since $X$ is a symplectic variety, $X$ is Cohen-Macaulay (cf. (1.1)). 
Moreover, $X$ is affine and 
$X \setminus U$ has codimension $\geq 4$ in $X$. 
Thus, by a depth argument, we see that 
$H^1(\mathcal{O}_U) = H^2(\mathcal{O}_U) = 0$. 
On the other hand,   
by Grothendieck's theorem [Gr]\footnote{The V-manifold case 
is reduced to the smooth case as follows. Let $W$ be an algebraic 
variety with quotient singularities (V-manifold). One can cover 
$W$ by finite affine open subsets $U_i$, $0 \leq i \leq n$ so 
that each $U_i$ admits an etale Galois cover $U'_i$ such 
that $U'_i = V_i/G_i$ with a smooth variety $V_i$ and a 
finite group $G_i$. It can be checked that, for each 
intersection $U_{i_0, ..., i_p} := U_{i_0} \cap ... \cap U_{i_p}$, 
Grothendieck's theorem holds. Now one has  
Grothendieck's theorem for $W$ by comparing two spectral 
sequences   
$$ E_1^{p,q} := \oplus_{i_0 < ... < i_p} H^q(U_{i_0, ..., i_p}, 
\tilde{\Omega}^{\cdot}_{U_{i_0, ..., i_p}}) \Longrightarrow 
H^{p+q}(W, \tilde{\Omega}^{\cdot}_W) $$ and 
$$ {E'}_1^{p,q} :=   \oplus_{i_0 < ... < i_p} H^q(U^{an}_{i_0, ..., i_p}, 
\mathbf{C}) \Longrightarrow 
H^{p+q}(W^{an}, \mathbf{C}).$$ } 
for $V$-manifolds, we 
have $\mathbf{H}^2(\tilde{\Omega}^{\cdot}_U) 
\cong \mathbf{H}^2(U^{an}, \mathbf{C})$. 
Now the result follows from the exact sequence above.  Q.E.D. 
\vspace{0.2cm}

\section{Prorepresentability of the Poisson deformation functors} 

Let $(X, \{\;,\;\})$ be a Poisson scheme. 
In this section, we shall prove that, in many important 
cases, $\mathrm{PD}_{(X, \{\;, \;\})}$ has a prorepresentable hull $R_X$ (cf. 
[Sch]), and it 
is actually prorepresentable, i.e.  $\mathrm{Hom}(R_X, \cdot) \cong 
\mathrm{PD}_{(X, \{\;, \;\})}(\cdot)$.  
Let $\mathcal{X}$ be a Poisson scheme over a local Artinian base $T$ and 
let $X$ be the central closed fiber. 
Let $G_{\mathcal{X}/T}$ be the sheaf of automorphisms of $\mathcal{X}/T$. 
More exactly, it is a sheaf on $X$ which associates to each open set $U \subset X$,  
the set of the automorphisms of the usual scheme $\mathcal{X}\vert_U$ over $T$ 
which induce the identity map on the central fiber $U = X\vert_U$. 
Moreover, let $PG_{\mathcal{X}/T}$ be the sheaf of {\em Poisson automorphisms} 
of  $\mathcal{X}/T$ as a subsheaf of $G_{\mathcal{X}/T}$.  
In order to show that $\mathrm{PD}_{(X, \{\;, \;\})}$ is prorepresentable, it 
is enough to prove that $H^0(X, PG_{\mathcal{X}/T}) \to H^0(X, PG_{\bar{\mathcal{X}}/\bar{T}})$ 
is surjective for any closed subscheme $\bar{T} \subset T$ and $\bar{\mathcal{X}} := 
\mathcal{X} \times_T \bar{T}$.         
Assume that $\mathcal{X}$ is smooth over $T$.  
We denote by  $\Theta_{\mathcal{X}/T}$ the relative tangent sheaf for $\mathcal{X} \to T$.   
Consider the Lichnerowicz-Poisson complex (cf. [Na 2, Section 2])  
$$ 0 \to \Theta_{\mathcal{X}/T} \stackrel{\delta_1}\to 
\wedge^2\Theta_{\mathcal{X}/T} \stackrel{\delta_2}\to 
\wedge^3\Theta_{\mathcal{X}/T} ... $$ and define 
$P\Theta_{\mathcal{X}/T} := \mathrm{Ker}(\delta_1)$. 
We denote by $\Theta^0_{\mathcal{X}/T}$ (resp. $P\Theta^0_{\mathcal{X}/T}$) 
the subsheaf of $\Theta_{\mathcal{X}/T}$ (resp. $P\Theta^0_{\mathcal{X}/T}$) 
which consists of the sections vanishing on the central closed fiber.   
\vspace{0.2cm}

{\bf Proposition (2.1)}(Wavrik): {\em There is an isomorphism of 
sheaves of sets  
$$  \alpha : \Theta^0_{\mathcal{X}/T} \cong G_{\mathcal{X}/T}. $$ 
Moreover, $\alpha$ induces an injection 
$$ P\Theta^0_{\mathcal{X}/T} \to PG_{\mathcal{X}/T}. $$}  

{\em Proof}.  Each local section $\varphi$ of 
$\Theta^0_{\mathcal{X}/T} $ is regarded as a derivation 
of $\mathcal{O}_{\mathcal{X}}$. Then we put 
$$ \alpha(\varphi) := id + \varphi + 1/2! (\varphi \circ \varphi) 
+ 1/3! (\varphi \circ \varphi \circ \varphi) + ... $$ 
By using the property 
$$ \varphi(fg) = f\varphi(g) + \varphi(f)g, $$ 
one can check that $\alpha(\varphi)$ is an automorphism 
of $\mathcal{X}/T$ inducing the identity map on the 
central fiber.  If $\varphi$ is a local section of 
$P\Theta^0_{\mathcal{X}/T}$, then $\varphi$ 
satisfies 
$$ \varphi(\{f,g\}) = \{f, \varphi(g)\} + \{\varphi(f), g\}.$$ 
By this property, one sees that $\alpha(\varphi)$ becomes 
a Poisson automorphism of $\mathcal{X}/T$.  
For the bijectivity of $\alpha$, see [Wav].    
\vspace{0.2cm}

{\bf Proposition (2.2)}. {\em In Proposition (2.1), if $\mathcal{X}$ 
is a Poisson deformation of a smooth symplectic variety $(X, \omega)$, then 
$\alpha$ induces an isomorphism} 
$$ P\Theta^0_{\mathcal{X}/T} \cong PG_{\mathcal{X}/T}. $$

{\em Proof}.  We only have to prove that the map is 
surjective. 
We may assume that $X$ is affine. Let 
$S$ be the Artinian local ring with $T = \mathrm{Spec}(S)$ 
and let $m$ be the maximal ideal of $S$. Put $T_n := 
\mathrm{Spec}(S/m^{n+1})$. The sequence 
$$T_0 \subset T_1 \subset ... \subset T_k $$ terminates 
at some $k$ and $T_k = T$.  We put $X_n := \mathcal{X} 
\times_T T_n$. Let $\phi$ be a section of $PG_{\mathcal{X}/T}$. 
One can write  
$$ \phi \vert_{X_1} = id + \varphi_1$$ with 
$\varphi_1 \in m\cdot P\Theta_X$.  
By the next lemma, 
$\varphi_1$ lifts to some 
$\tilde{\varphi}_1 \in P\Theta_{\mathcal{X}/T}$. 
Then one can write 
$$\phi\vert_{X_2} = \alpha(\tilde{\varphi}_1)\vert_{X_2} 
+ \varphi_2$$ with $\varphi_2 \in m^2\cdot P\Theta_X$. 
Again, by the lemma, $\varphi_2$ lifts to some 
$\tilde{\varphi}_2 \in P\Theta_{\mathcal{X}/T}$. 
Continue this operation and we finally conclude that 
$$ \phi = \alpha(\tilde{\varphi}_1 + \tilde{\varphi}_2 + ... ). $$ 

{\bf Lemma (2.3)}. {\em Let $\mathcal{X} \to T$ be a Poisson 
deformation of a smooth symplectic variety $(X, \omega)$ over a local 
Artinian base $T = \mathrm{Spec}(S)$. Let $\bar{T} \subset T$ be a closed 
subscheme and put $\bar{\mathcal{X}} := \mathcal{X} \times_T \bar{T}$.  
Then the restriction map 
$$P\Theta_{\mathcal{X}/T} \to P\Theta_{\bar{\mathcal{X}}/\bar{T}} $$ 
is surjective.} 
\vspace{0.2cm}

{\em Proof}. We may assume that $X$ is affine.  The Lichnerowicz-Poisson complex 
$(\wedge^{\geq 1}\Theta_{\mathcal{X}/T}, \delta)$ is 
identified with the truncated de Rham complex 
$(\Omega^{\geq 1}_{\mathcal{X}/T}, d)$ by the symplectic 
2-form $\omega$ (cf. [Na 2], Section 2). 
There is a distinguished triangle 
$$ \Omega^{\geq 1}_{\mathcal{X}/T} \to \Omega^{\cdot}_{\mathcal{X}/T} 
\to \mathcal{O}_{\mathcal{X}} \to 
\Omega^{\geq 1}_{\mathcal{X}/T}[1],$$ and it induces   
an exact sequence  
$$ ... \to HP^i(\mathcal{X}/T) \to H^i(X^{an}, S) \to 
H^i(X, \mathcal{O}_{\mathcal{X}}) \to ... $$ 
In particular, we have an exaxt sequence 
$$ 0 \to K \to HP^1(\mathcal{X}/T) \to H^1(X^{an}, S) \to 0,$$ 
where $$K: = \mathrm{Coker}[H^0(X^{an}, S) \to H^0(X, \mathcal{O}_{\mathcal{X}})].$$  
Similarly for $\bar{\mathcal{X}}$, we have an exact sequence 
$$ 0 \to \bar{K} \to HP^1(\bar{\mathcal{X}}/\bar{T}) \to 
H^1(X^{an}, \bar{S}) \to 0$$ 
with $$\bar{K} := \mathrm{Coker}[H^0(X^{an}, \bar{S}) \to 
H^0(X, \mathcal{O}_{\bar{\mathcal{X}}})].$$
Since the restriction maps $K \to \bar{K}$ and 
$H^0(X^{an}, S) \to H^0(X^{an}, \bar{S})$ are both surjective, 
the restriction map $HP^1(\mathcal{X}/T) \to HP^1(\bar{\mathcal{X}}/\bar{T})$ 
is surjective. 
Finally, note that 
$HP^1(\mathcal{X}/T) = H^0(X, P\Theta_{\mathcal{X}/T})$ 
and 
$HP^1(\bar{\mathcal{X}}/\bar{T}) = H^0(X, P\Theta_{\bar{\mathcal{X}}/\bar{T}})$.
\vspace{0.2cm}

{\bf Proposition (2.4)}. {\em In the same assumption in 
Lemma (2.3),  if the restriction map 
$$ H^0(X, P\Theta_{\mathcal{X}/T}) \to 
H^0(X, P\Theta_{{\bar{\mathcal{X}}}/\bar{T}})$$ 
is surjective, then the restriction map 
$$ H^0(X, PG_{{\mathcal{X}}/T}) \to 
H^0(X, PG_{{\bar{\mathcal{X}}}/\bar{T}})$$ 
is surjective. }  
\vspace{0.2cm}

{\em Proof}. If the map 
$$ H^0(X, P\Theta_{\mathcal{X}/T}) \to 
H^0(X, P\Theta_{{\bar{\mathcal{X}}}/\bar{T}})$$ 
is surjective, 
$$ H^0(X, P\Theta^0_{\mathcal{X}/T}) \to 
H^0(X, P\Theta^0_{{\bar{\mathcal{X}}}/\bar{T}})$$ 
is surjective. Then the result follows from Proposition 
(2.2). 
\vspace{0.2cm}

{\bf Corollary (2.5)}. {\em The Poisson deformation 
functor $PD_{(X, \{\;, \;\})}$ for a symplectic variety 
$(X, \omega)$, is prorepresentable 
in the following two cases:} 

(1) {\em $X$ is convex (i.e. $X$ has a birational projective morphism 
to an affine variety), and admits only terminal singularities.} 

(2) {\em $X$ is affine, and $H^1(X^{an}, \mathbf{C}) = 0$.} 
\vspace{0.2cm}

{\em Proof}. First, we must show that 
$\dim \mathrm{PD}_{(X, \{\;,\;\})}(\mathbf{C}[\epsilon]) < \infty$. Let $U$ be the smooth 
part of $X$. In the case (1), we have 
$\mathrm{PD}_{(X, \{\;,\;\})}(\mathbf{C}[\epsilon]) = H^2(U^{an}, \mathbf{C})$; hence  
$\mathrm{PD}_{(X, \{\;,\;\})}(\mathbf{C}[\epsilon])$ is a finite dimensional $\mathbf{C}$-vector 
space. For the case (2), the finiteness is proved in Proposition (1.10).  
Assume that $\mathcal{X} \to T$ is a Poisson deformation 
of $X$ with a local Artinian base. Let $\bar{T}$ be a closed 
subscheme of $T$ and let $\bar{\mathcal{X}} \to \bar{T}$ be 
the induced Poisson deformation of $X$ over $\bar{T}$.  
Let $\mathcal{U} \subset \mathcal{X}$ (resp. $\bar{\mathcal{U}} \subset 
\bar{\mathcal{X}}$) be the open locus where the map 
$\mathcal{X} \to T$ (resp. $\bar{\mathcal{X}} \to \bar{T}$) is 
smooth. Let $j$ be the inclusion map of 
$\mathcal{U}$ to $\mathcal{X}$. Since $j_*\mathcal{O}_{\mathcal{U}} 
= \mathcal{O}_{\mathcal{X}}$, a Poisson automorphism of $\mathcal{U}$ 
(which induces the identity on the closed fiber) uniquely extends to 
that of $\mathcal{X}$. Therefore, we have an isomorphism 
$$ H^0(\mathcal{X}, PG_{\mathcal{X}/T}) \cong 
H^0(\mathcal{U}, PG_{\mathcal{U}/T}). $$ 
Similarly, we have 
$$ H^0(\bar{\mathcal{X}}, PG_{\bar{\mathcal{X}}/\bar{T}}) \cong 
H^0(\bar{\mathcal{U}}, PG_{\bar{\mathcal{U}}/\bar{T}}).$$ 
By Proposition (2.4), it suffices to show that the 
restriction map 
$$ H^0(U, P\Theta_{\mathcal{U}/T}) \to 
H^0(U, P\Theta_{\bar{\mathcal{U}}/\bar{T}})$$   
is surjective. 

For the case (1), we have already proved the surjectivity 
in [Na 2], Theorem 14. 
Let us consider the case (2). 
Note that $H^0(U, P\Theta_{\mathcal{U}/T}) 
\cong \mathbf{H}^1(U, \Theta^{\geq 1}_{\mathcal{U}/T})$, 
where $(\Theta^{\geq 1}_{\mathcal{U}/T}, \delta)$ is 
the Lichnerowicz-Poisson complex for $\mathcal{U}/T$. 
As in the proof of Lemma (2.3), the Lichnerowicz-Poisson 
complex is identified with the truncated de Rham complex 
$(\Omega^{\geq 1}_{\mathcal{U}/T}, d)$, 
and it induces the exact sequence 
$$ 0 \to K \to \mathbf{H}^1(U, \Omega^{\geq 1}_{\mathcal{U}/T}) 
\to H^1(U^{an}, S), $$ 
where $S$ is the affine ring of $T$, and 
$K := \mathrm{Coker}[H^0(U^{an}, S) \to H^0(U, \mathcal{O}_{\mathcal{U}})]$. 
We shall prove that $H^1(U^{an}, S) = 0$. 
Since $H^1(U^{an}, S) = H^1(U^{an}, \mathbf{C}) \otimes S$, 
it suffices to show that $H^1(U^{an}, \mathbf{C}) = 0$. 
Let $f: \tilde{X} \to X$ be a resolution of $X$ such that 
$f^{-1}(U) \cong U$ and the  
exceptional locus $E$ of $f$ is a divisor with only simple normal crossing. 
One has the exact sequence 
$$ H^1(\tilde{X}^{an}, \mathbf{C}) \to H^1(U^{an}, \mathbf{C}) 
\to H^2_E(\tilde{X}^{an}, \mathbf{C}) \to H^2(\tilde{X}^{an}, \mathbf{C}), $$ 
where the first term is zero because $X$ has only rational singularities 
and $H^1(X^{an}, \mathbf{C}) = 0$. We have to prove that 
$H^2_E(\tilde{X}^{an}, \mathbf{C}) \to H^2(\tilde{X}^{an}, \mathbf{C})$ 
is an injection. Put $n := \dim X$; then, $H^2_E(\tilde{X}^{an}, \mathbf{C})$ 
is dual to the cohomology $H^{2n-2}_c(E^{an}, \mathbf{C})$ with 
compact support (cf. the proof of Proposition 2 of [Na 3]). 
Let $E = \cup E_i$ be the irreducible decomposition of $E$. The 
$p$-multiple locus of $E$ is, by definition, the locus of points 
of $E$ which are contained in the intersection of some $p$ different 
irreducible components of $E$. Let $E^{[p]}$ be the normalization 
of the $p$-multiple locus of $E$. For example, $E^{[1]}$ is the 
disjoint union of $E_i$'s, and $E^{[2]}$ is the normalization 
of the singular locus of $E$. 
There is an exact sequence 
$$ 0 \to \mathbf{C}_{E} \to \mathbf{C}_{E^{[1]}} \to 
\mathbf{C}_{E^{[2]}} \to ... $$
By using this exact sequence, we see that $H^{2n-2}_c(E^{an}, \mathbf{C})$ 
is a $\mathbf{C}$-vector space whose dimension equals the 
number of irreducble components of $E$. 
By the duality, we have 
$$H^2_E(\tilde{X}^{an}, \mathbf{C}) = \oplus \mathbf{C}[E_i]$$ 
and the map $H^2_E(\tilde{X}^{an}, \mathbf{C}) \to 
H^2(\tilde{X}^{an}, \mathbf{C})$ is an injection. 
Therefore, $H^1(U^{an}, \mathbf{C}) = 0$. 
We now know that $$H^0(U, P\Theta_{\mathcal{U}/T}) \cong 
K.$$ Similarly, we have $$H^0(U, P\Theta_{\bar{\mathcal{U}}/\bar{T}}) 
\cong \bar{K},$$ where $\bar{K} := 
\mathrm{Coker}[H^0(U, \bar{S}) \to H^0(U, \mathcal{O}_{\bar{\mathcal{U}}})]$ 
and $\bar{S}$ is the affine ring of $\bar{T}$. 
Since the restriction maps $H^0(X, \mathcal{O}_{\mathcal{X}}) \to 
H^0(U, \mathcal{O}_{\mathcal{U}})$ and 
$H^0(X, \mathcal{O}_{\bar{\mathcal{X}}}) \to 
H^0(U, \mathcal{O}_{\bar{\mathcal{U}}})$ are both isomorphisms, 
the restriction map 
$H^0(U, \mathcal{O}_{\mathcal{U}}) \to 
H^0(U, \mathcal{O}_{\bar{\mathcal{U}}})$ is surjective; hence 
the map $K \to \bar{K}$ is also surjective. Q.E.D. 
\vspace{0.15cm}

{\bf Remark (2.6)}. The results in this section equally hold in the 
complex analytic category. For example, let $(X,p)$ be the germ of a 
symplectic variety $X$ at $p \in X$, and let $f: (Y, E) \to (X,p)$ be a 
crepant, projective partial resolution of $(X,p)$ where $E = f^{-1}(p)$. 
Assume that $Y$ has only terminal singularities. 
Then (2.5) holds for $(X,p)$ and $(Y,E)$.   

\section{Symplectic automorphism and universal Poisson deformations} 

Let $S$ be the same as in (1.2), and put  
$V := (S,0) \times (\mathbf{C}^{2n-2},0)$. By the symplectic 
2-form $\omega := (p_1)^*\omega_S + (p_2)^*\omega_{\mathbf{C}^{2n-2}}$, 
the germ $V$ becomes a symplectic variety. Let $(\tilde{S}, F) 
\to (S,0)$ be the minimal resolution and put 
$\tilde{V} := (\tilde{S}, F) \times (\mathbf{C}^{2n-2},0)$.  
In this section, we construct explicitly the universal Poisson 
deformations of $V$ and $\tilde{V}$, and study the natural action 
on them induced by a symplectic automorphism of $V$. 
Let $\g$ be the complex simple Lie algebra of the same type as $S$. Fix a Cartan subalgebra 
$\h$ of $\g$ and consider the adjoint quotient map $\g \to \h/W$, where $W$ is the Weyl group of  
$\g$. By [Slo], a transversal slice  ${\mathcal S}$ of  $\g$ at the sub-regular nilpotent orbit 
gives the semi-universal flat deformation $\mathcal{S} \to \h/W$ of $S$ (at $0 \in \h/W$).
Let $\g_{reg}$ be the open set of $\g$ where this map is smooth. Then $\g_{reg} \to \h/W$ admits a 
relative symplectic 2-form called the Kostant-Kirillov 2-form.   
Let $\mathcal{S}_{reg}$ 
be the open subset of $\mathcal{S}$ where the map $\mathcal{S} \to \h/W$ is smooth.  
The Kostant-Kirillov 2-form on $\g_{reg}$ restricts to give a relative symplectic 2-form on 
$\mathcal{S}_{reg}$ and  
makes the map $\mathcal{S} \to \h/W$ a Poisson deformation of $S$. 

On the other hand, the base change $\g \times_{\h/W}\h \to \h$ has a simultaneous resolution 
$$\mu: G \times^B \b \to \g \times_{\h/W}\h,$$ 
where $G$ is the adjoint group of $\g$ and $B$ is a Borel 
subgroup of $G$ such that $\h \subset \b$ (cf. [Slo]). 
The pullback of the Kostant-Kirillov 2-form 
gives a relative symplectic 2-form 
$\omega_F \in \Gamma(G \times^B \b, \Omega^2_{G \times^B \b/\h})$.   
If we put $\tilde{\mathcal S} := \mu^{-1}(\mathcal{S} \times_{\h/W} \h)$, then 
$$\mu\vert_{\tilde{\mathcal S}}: \tilde{\mathcal S} \to \mathcal{S} \times_{\h/W}\h$$ 
is a simultaneous resolution of $\mathcal{S} \times_{\h/W}\h \to \h$. 
Let $f$ be the composite of two maps $\tilde{\mathcal S} \to \mathcal{S} \times_{\h/W}\h$ 
and $\mathcal{S} \times_{\h/W}\h \to \h$. Then $\omega_f := \omega_F\vert_{\tilde{\mathcal S}}$ 
gives a relative symplectic 2-form for $f$ (cf. [Ya]). 
\vspace{0.2cm}

{\bf Proposition (3.1)} {\em (1) The universal Poisson deformations of 
$S$ and $\tilde{S}$ are respectively given by ${\mathcal S} \to \h/W$ and 
$\tilde{\mathcal S} \to \h$.} 

{\em (2) The universal Poisson deformations of $V$ and $\tilde{V}$ are 
respectively given by ${\mathcal S} \times ({\mathbf C}^{2n-2},0) \to \h/W$ 
and $\tilde{\mathcal S} \times ({\mathbf C}^{2n-2},0) \to \h$.} 

{\em Proof}. The Poisson deformation $\mathcal{S} \to \h/W$ 
is universal at $0 \in \h/W$. In fact, 
there is an exact sequence (cf. the 
latter part of \S 1 after (1.8)) 
$$ 0 \to \mathrm{PD}_{lt, S}(\mathbf{C}[\epsilon]) \to \mathrm{PD}_S(\mathbf{C}[\epsilon]) 
\to T^1_S \to 0.$$ For the definitions of $\mathrm{PD}$ and $\mathrm{PD}_{lt}$, see (1.1) and (1.8).  
By Proposition (1.11), we have $\mathrm{PD}_{lt, S}(\mathbf{C}[\epsilon]) \cong 
H^2(S, \mathbf{C}) = 0$. The map $\mathrm{PD}_S(\mathbf{C}[\epsilon]) \to T^1_S$ is an isomorphism. 
Since $\mathcal{S} \to \h/W$ is a semi-universal flat deformation of $S$, the Kodaira-Spencer map 
$T_{\h/W,0} \to T^1_S$ is an isomorphism. The Kodaira-Spencer map factorizes as 
$T_{\h/W,0} \to \mathrm{PD}_S(\mathbf{C}[\epsilon]) \to T^1_S$; hence the Poisson Kodaira-Spencer map 
$T_{\h/W,0} \to \mathrm{PD}_S(\mathbf{C}[\epsilon])$ is an isomorphism. This fact together with (2.6) 
implies the universality of the Poisson deformation.     
Now let us consider the map $\tilde{\mathcal S} \to \h$. 
By [Slo], it is semi-universal as a usual flat deformation of $\tilde{S}$. Therefore, the Kodaira-Spencer 
map $T_{\h,0} \to H^1(\tilde{S}, \Theta_{\tilde S})$ is an isomorphism. Moreover, this map factorizes 
as $T_{\h,0} \to H^2(\tilde{S}, \mathbf{C}) \to H^1(\tilde{S}, \Theta_{\tilde S})$, where 
the map $T_{\h,0} \to H^2(\tilde{S}, \mathbf{C})$ is the Poisson Kodaira-Spencer map. By the 
symplectic 2-form, $\Theta_{\tilde S}$ and $\Omega^1_{\tilde S}$ are identified. Then,  
the map $H^2(\tilde{S}, \mathbf{C}) \to H^1(\tilde{S}, \Theta_{\tilde S})$ coincides with the natural  
isomorphism $H^2(\tilde{S}, \mathbf{C}) \to H^1(\tilde{S}, \Omega^1_{\tilde S})$. 
Therefore, the Poisson Kodaira-Spencer map 
$T_{\h,0} \to H^2(\tilde{S}, \mathbf{C})$ is an isomorphism. This fact together with (2.6) implies that 
$f: \tilde{\mathcal S} \to \h$ is the universal Poisson deformation of $\tilde{S}$. 
Let us now consider the Poisson deformations of $\tilde{V}$. 
The tangent space $\mathrm{PD}_{\tilde V}({\mathbf C}[\epsilon])$ of the 
Poisson deformation functor is isomorphic to 
$H^2(\tilde{S} \times {\mathbf C}^{2n-2}, {\mathbf C}) = H^2({\tilde S}, {\mathbf C})$.     
Since $\mathrm{PD}_{(\tilde S, F)}({\mathbf C}[\epsilon]) \cong H^2(\tilde{S}, {\mathbf C})$, this means 
that $$\tilde{{\mathcal S}} \times {\mathbf C}^{2n-2} \stackrel{f \circ p_1}\to \h$$ is the universal Poisson 
deformation of $\tilde{V}$ at $0 \in \h$. 
Moreover, the map $$\mathcal{S} \times {\mathbf C}^{2n-2} \to \h/W$$ is the universal Poisson 
deformation of $V$ at $0 \in \h/W$. In fact, the map 
$\mathcal{S} \to \h/W$ is the universal Poisson deformation of $S$. By Lemma (1.6), any 1-st 
order Poisson deformation is the product of a 1-st order Poisson deformation of $S$ and 
$(\mathbf{C}^{2n-2},0)$. Then, the Poisson Kodaira-Spencer map 
$T_{\h/W,0} \to \mathrm{PD}_{V}(\mathbf{C}[\epsilon])$ is an isomorphism. Q.E.D. 
\vspace{0.2cm}

Let $$i: V \to V$$ be a symplectic automorphism of $V$. 
The map $i$ lifts to a symplectic automorphism 
$$ \tilde{i}: \tilde{V} \to \tilde{V}$$ so that the following diagram commutes 
\begin{equation} 
\begin{CD}  
\tilde{V} @>{\tilde i}>> \tilde{V} \\ 
@VVV @VVV \\ 
V @>{i}>> V.   
\end{CD} 
\end{equation} 
Correspondingly, we have a commutative diagram of functors:  

\begin{equation} 
\begin{CD}  
\mathrm{PD}_{\tilde V} @>{\tilde{i}_*}>> 
\mathrm{PD}_{\tilde V} \\ 
@VVV @VVV \\ 
\mathrm{PD}_{V} @>{i}>> \mathrm{PD}_{V}.    
\end{CD} 
\end{equation} 

By the (formal) universality of $\mathrm{PD}_{V}$ and 
$\mathrm{PD}_{\tilde V}$ (cf.(2.5), (2.6)), we have a commutative diagram 
  
\begin{equation} 
\begin{CD} 
\hat{\h} @>{\tilde{\iota}}>> \hat{\h}\\ 
@VVV @VVV \\ 
\hat{\h/W} @>{\iota}>> \hat{\h/W},     
\end{CD} 
\end{equation} 

\noindent where $\hat{\h}$ and $\hat{\h/W}$ are the 
formal completions of $\h$ and $\h/W$ at the origins. 
\vspace{0.2cm}

{\bf Proposition (3.2)} {\em The quotient space $\h/W$ has a linear 
structure so that the commutative diagram above is obtained from a 
commutative diagram of linear spaces} 

\begin{equation} 
\begin{CD} 
\h @>>> \h\\ 
@VVV @VVV \\ 
\h/W @>>> \h/W    
\end{CD} 
\end{equation} 
 
\noindent {\em where both horizontal maps are linear maps. Moreover, 
the horizontal map $\h \to \h$ is induced by a graph 
automorphism of the Dynkin diagram of $\g$.} 
\vspace{0.2cm}

{\em Proof}. Let us consider the Poisson deformation 
$\tilde{\mathcal S} \times ({\mathbf C}^{2n-2},0) \to \h$. 
The relative symplectic 2-form $\omega_f + \omega_{{\mathbf C}^{2n-2}}$ 
defines a 2-nd cohomology class of each fiber 
$\tilde{\mathcal S}_t \times ({\mathbf C}^{2n-2},0)$, $t \in \h$. 
Since $H^2(\tilde{S}_t \times \mathbf{C}^{2n-2}, \mathbf{C})$ 
is identified with $H^2(\tilde{V}, \mathbf{C})$, 
one can define a period map (cf. [G-K], [Ya])  
$$ p: \h \to H^2(\tilde{V}, \mathbf{C}) \cong H^2(\tilde{S}, \mathbf{C}).$$ 
Similarly one can define a period map 
$$ p_B: \h \to H^2(T^*(G/B), \mathbf{C})$$ for 
the Poisson deformation $F: (G \times^B \b)\times (\mathbf{C}^{2n-2},0) \to \h$ 
by using the relative symplectic 2-form $\omega_F + \omega_{{\mathbf C}^{2n-2}}$. 
Since $\omega_f = \omega_F\vert_{\tilde{\mathcal S}}$, 
the period map $p$ is the composite of $p_B$ and the natural 
restriction map $H^2(T^*(G/B), \mathbf{C}) \to H^2(\tilde{S}, \mathbf{C})$. 
This restriction map is an isomorphism since $\g$ is simply-laced. 
Note that $W$ has monodromy actions on $H^2(T^*(G/B), \mathbf{C})$ and $H^2(\tilde{S}, \mathbf{C})$ 
([Slo 2], 4.2, 4.3, 4.4). 
By [Ya, Section 3] the period map $p_B$ is a $W$-equivariant linear isomorphism; hence 
$p$ is also a $W$-equivariant linear isomorphism. 
The description of $p_B$ is as follows. First of all, 
the nilpotent cone $N$ of $\g$ is resolved by the 
Springer map $\mu_0: T^*(G/B) \to N$. The transversal slice $S$ is contained 
in $N$ and $\tilde{S} = \mu_0^{-1}(S)$. 
There is an isomorphism 
$$ \h^* \to H^2(T^*(G/B), \mathbf{C}).$$ 
The construction is as follows. Let $H \subset B$ be the maximal torus corresponding 
to $\h$. Then there is a canonical isomorphism (cf. [Na 5, (P3)]) 
$$ \mathrm{Hom}_{alg.gp}(H, \mathbf{C}^*)\otimes \mathbf{C} \cong 
\mathrm{Pic}(G/B) \otimes \mathbf{C}.$$ 
The left hand side is $\h^*$ and right hand side is isomorphic to 
$H^2(G/B, \mathbf{C})$. Since $H^2(G/B, \mathbf{C}) \cong 
H^2(T^*(G/B), \mathbf{C})$, we have an isomorphism $\h^* \to H^2(T^*(G/B), \mathbf{C})$. 
The Cartan subalgebra $\h$ is identified with its dual 
$\h^*$ by the Killing form of $\g$. 
By \S 3 of [Ya] the period map $p_B$ coincides with the composite of two maps: 
$$\h \to \h^* \to H^2(T^*(G/B), \mathbf{C}).$$ 
The automorphism $\tilde{i}$ of $\tilde{V}$ induces an isomorphism 
$${\tilde i}^*: H^2(\tilde{V}, \mathbf{C}) \to H^2(\tilde{V}, \mathbf{C}).$$ 
By the identification $H^2(\tilde{V}, \mathbf{C}) \cong H^2(\tilde{S}, \mathbf{C})$, 
the map ${\tilde i}^*$ is regarded as an automorphism of $H^2(\tilde{S}, \mathbf{C})$. 
By the definition of $\tilde{\iota}$ we have a commutative diagram 
 
\begin{equation} 
\begin{CD}  
\hat{\h} @>>> \h @>{p}>> H^2(\tilde{S}, \mathbf{C}) \\ 
@A{\tilde{\iota}}AA  @.   @V{{\tilde i}^*}VV \\ 
\hat{\h} @>>> \h @>{p}>> H^2(\tilde{S}, \mathbf{C}).   
\end{CD} 
\end{equation} 

Define a linear map $\tilde{\iota}_{\h}: \h \to \h$ by $p^{-1} \circ (\tilde{i}^*)^{-1} 
\circ p$. Then we have a commutative diagram 

\begin{equation} 
\begin{CD}  
\hat{\h} @>>> \h \\ 
@V{\tilde{\iota}}VV @V{\tilde{\iota}_{\h}}VV \\ 
\hat{\h} @>>> \h  
\end{CD} 
\end{equation} 
     
We shall prove that $\tilde{\iota}_{\h}$ is induced by a graph automorphism of $\g$. 
Let $\Phi \subset \h$ be the (co)root system for $\g$. The choice of $B$ determines a 
base $\Delta$ of $\Phi$. Define 
$$\Gamma := \{\phi \in \mathrm{Aut}(\Phi); \phi(\Delta) = \Delta\}.$$ 
Let $C_i$ be a $(-2)$-curve on $\tilde{S}$ and let $[C_i] \in H^2(\tilde{S}, \mathbf{C})$ 
be its class. Define 
$$ \Phi' := \{C := \Sigma a_i[C_i]; a_i \in \mathbf{Z}, \; C^2 = -2\}.$$ 
Then $\Phi'$ becomes a root system and $\Delta' := \{[C_i]\}$ forms a 
base of $\Phi'$. Define 
$$\Gamma' := \{\phi \in \mathrm{Aut}(\Phi'); \phi(\Delta') = \Delta' \}.$$ 
The period map $p$ sends $\Delta$ to $\Delta'$ up to a 
non-zero constant. Since ${\tilde i}^* \in \Gamma'$, we have $\tilde{\iota}_{\h} \in \Gamma$. 
The Weyl group $W$ of $\g$ is a normal subgroup of $\mathrm{Aut}(\Phi)$ and 
$\mathrm{Aut}(\Phi)$ is the semi-direct product 
of $W$ and $\Gamma$.
This means that $\tilde{\iota}_{\h}$ descends to an automorphism ${\iota}_{\h/W}$ of 
$\h/W$. Since $W$ is a finite reflection group, $\h/W$ is an affine 
space.  By [Slo, 8.8, Lemma 1], one can choose a linear structure of 
$\h/W$ so that ${\iota}_{\h/W}$ is a linear map. 
\vspace{0.2cm}

\section{Global sections of the local system}  

(4.1) {\bf Monodromy of $R^2\pi^{an}_*\mathbf{C}$}

As in (1.2)-(1.5), we shall consider a symplectic 
variety $(X,\omega)$ whose singularities are locally isomorphic to 
$(S,0) \times ({\mathbf C}^{2n-2},0)$.  We use the same notation in 
section 1. 
Let $\pi: Y \to X$ be the minimal resolution. By definition, 
$\pi^{an}$ is locally a product of the minimal resolution $\tilde{S} \to S$ 
and the $2n-2$ dimensional disc $\Delta^{2n-2}$. If $S$ is of type 
$A_r$, $D_r$ or $E_r$, then, for each $p \in \Sigma$, the fiber $(\pi^{an})^{-1}(p)$ 
has $r$ irreducible components and each of them is isomorphic 
to $\mathbf{P}^1$.  Let $E$ be the $\pi$-exceptional locus 
and let $m$ be the number of irreducible components of $E$. 
We have $m \le r$; but $m \ne r$ in general.  The local system 
$R^2\pi^{an}_*\mathbf{C}$ on $\Sigma$ may 
possibly have monodromies. Let $\gamma$ be a closed loop in $\Sigma$ starting from 
$p \in \Sigma$. Then we have a monodromy transformation along $\gamma$: 
$$H^2((\pi^{an})^{-1}(p), \mathbf{C}) \to H^2((\pi^{an})^{-1}(p), \mathbf{C}).$$ 
Since $H^2((\pi^{an})^{-1}(p), \mathbf{C}) \cong H^2(\tilde{S}, \mathbf{C})$, the 
monodromy transformation is an automorphism of $H^2(\tilde{S}, \mathbf{C})$. 
Let $F$ be the exceptional divisor of the minimal resolution $\tilde{S} \to S$ 
and let $F = \cup F_i$ be the irreducible decomposition. Then $\{[F_i]\}$ is a  
basis of $H^2(\tilde{S}, \mathbf{C})$. The monodromy transformation permutes  
$[F_i]$'s without changing the intersection numbers. Therefore, the monodromy 
transformation comes from a graph automorphism of the Dynkin diagram associated 
with $S$.  Let us observe the graph automorphisms of various Dynkin diagrams. 
In the $(A_r)$-case, the Dynkin 
diagram 

\begin{picture}(300, 20)(0, 0) 
\put(25, -10){1}\put(30, -3){$\circ$}\put(35, 0){\line(1, 0){25}} 
\put(65, -3.5){- - -}\put(90, 0){\line(1, 0){15}}
\put(105, -3){$\circ$}\put(110, 0){\line(1, 0){10}}
\put(125, -3.5){- - -}\put(150, 0)
{\line(1, 0){55}}\put(202, -10){r}\put(207, -3){$\circ$}    
\end{picture} 
\vspace{0.7cm}

\noindent has an automorphism $\sigma_1$ of order 2 which sends each $i$-th vertex to 
the $r+1-i$-th vertex.  Hence, there are two possibilities for $m$; namely,  
$$ m = r, \;\mathrm{or}\;  r- [r/2]. $$ 
The Dynkin diagram of type $D_r$ 
 
\begin{picture}(300, 20)(0, 0) 
\put(25, 3){1} \put(30, 10){$\circ$}\put(30, -13){$\circ$}\put(25, -20){2} 
\put(35, 10){\line(1, -1){10}}\put(35, -10){\line(1, 1){10}}
\put(47, -3){$\circ$}\put(55, 0){\line(1, 0){25}}
\put(85, -3.5){- - -}\put(110, 0){\line(1, 0){55}}\put(167, -3)
{$\circ$}\put(162, -10){r}
\end{picture} 
\vspace{1.0cm}

\noindent has an automorphism $\sigma_2$ of order 2, which sends the $1$-st vertex to 
the $2$-nd one. Especially when $r = 4$, it has another automorphism $\tau$ of order 3
which permutes mutually the $1$-st vertex, the $2$-nd one and $3$-rd one.  
Hence, in the $(D_4)$-case, there are three possibilities 
for $m$ 
$$  m= 4, 3 \;\mathrm{or}\; 2,$$
and, in the $(D_r)$-case with $r > 4$, there are two possibilities for $m$ 
$$  m = r \;\mathrm{or}\; r-1.$$ 

Finally, let us consider the $(E_6)$-case. 
  
\begin{picture}(300, 20)
\put(25, -10){1}\put(30, -3){$\circ$}\put(35, 0){\line(1, 0){20}}
\put(52, -10){2}\put(57, -3){$\circ$}\put(65, 0){\line(1, 0){20}}
\put(82,-10){3}
\put(87, -3){$\circ$}\put(90, -5){\line(0, -1){10}}
\put(82,-20){4}
\put(87, -20){$\circ$}\put(95, 0){\line(1, 0){20}}
\put(112, -10){5}\put(117, -3){$\circ$}\put(125, 0){\line(1, 0){20}}
\put(142, -10){6}\put(147, -3){$\circ$} 
\end{picture} 
\vspace{0.7cm}

The diagram has an automorphism $\sigma_3$ of order 2, which sends the 
$1$-st vertex to the $6$-th one and the $2$-nd one to the $5$-th 
one.  There are two possibilities for $m$ 
$$ m = 6, \;\mathrm{or}\; 4.$$ 
Since there are no symmetries for the diagrams of type $(E_7)$, $(E_8)$,  
we conclude that $m = r$ in these cases. 

Let $\gamma$ be a closed loop in $\Sigma$ starting from $p \in \Sigma$. 
In (1.7), we have chosen a sequence of points $p_i$ ($1 \le i \le k$) on $\gamma$ and   
have made a sequence of symplectic isomorphisms $(X^{an}, p_{i-1}) \cong (X^{an}, p_i)$. 
The composite of them finally defines a symplectic automorphism 
$$ i_{\gamma}: (X^{an}, p) \cong (X^{an}, p).$$  
Here we shall describe the monodromy transformation of $R^2\pi^{an}_*\mathbf{C}$ along $\gamma$ in 
terms of a symplectic automorphism of $(Y^{an}, \pi^{-1}(p))$.  
For each open set $V \subset X^{an}$, we associate the $\mathbf{C}$-vector space which consists of  
all 1-st order Poisson deformations of $(\pi^{an})^{-1}(V)$. The sheaf determined by 
this presheaf is isomorphic to $R^2\pi^{an}_*\mathbf{C}$ (cf. [Na 2]).    
The symplectic isomorphisms $(X^{an}, p_{i-1}) \cong (X^{an}, p_i)$ induce 
symplectic isomorphisms $(Y^{an}, (\pi^{an})^{-1}(p_{i-1})) \cong 
(Y^{an}, (\pi^{an})^{-1}(p_i))$ because $(Y^{an}, (\pi^{an})^{-1}(p_i))$ is a unique crepant resolution of 
$(X^{an}, p_i)$. The sequence of them finally defines a 
symplectic automorphism  
$$ \tilde{i}_{\gamma} : (Y^{an}, (\pi^{an})^{-1}(p) \cong (Y^{an}, (\pi^{an})^{-1}(p)). $$ 
Note that $\tilde{i}_{\gamma}$ is a (unique) lift of ${i}_{\gamma}$ to an automorphism 
of $(Y^{an}, (\pi^{an})^{-1}(p))$. 
The map $\tilde{i}_{\gamma}$ induces an automorphism of $(R^2\pi^{an}_*\mathbf{C})_p$, 
which is nothing but the monodromy transformation of $R^2\pi^{an}_*\mathbf{C}$ 
along $\gamma$. 
The identification $(X^{an}, p) \cong (S,0) \times ({\mathbf C}^{2n-2},0)$ naturally 
lifts to the identification of $(Y^{an}, (\pi^{an})^{-1}(p))$ with $(\tilde{S}, F) \times 
(\mathbf{C}^{2n-2},0)$. Then, $(R^2\pi^{an}_*\mathbf{C})_p$ can be identified with   
$H^2(\tilde{S}, \mathbf{C})$. 

The following is the main 
result in this section.
\vspace{0.2cm}

{\bf Proposition (4.2)}. {\em The following equality holds:} 
$$ \dim_{\mathbf{C}}H^0(\Sigma, \mathcal{H}) =  m.$$ 
\vspace{0.15cm}

{\em Proof}. Let $\gamma$ be a closed loop starting from $p \in \Sigma$. 
As in (1.7), we choose admissible covers $\{U_i\}$ of $\gamma$ and 
points $p_i \in \Gamma$. By (1.7) and (4.1), the monodromy transformations 
of $\mathcal{H}_p$ and $(R^2\pi^{an}_*\mathbf{C})_p$ along $\gamma$, are 
described in terms of symplectic automorphisms $${i}_{\gamma}: (X^{an}, p) 
\to (X^{an},p)$$ and $$\tilde{i}_{\gamma}: (Y^{an}, (\pi^{an})^{-1}(p)) \to 
(Y^{an}, (\pi^{an})^{-1}(p)).$$ Apply Proposition (3.2) to these symplectic automorphisms.   
Then the sheaf $R^2\pi^{an}_*\mathbf{C}$ is a local system of the {\bf C}-module $\h$, 
and $\mathcal{H}$ is a local system of the {\bf C}-module $\h/W$. Moreover, 
their monodromies along $\gamma$ are given by the horizontal maps $\h \to \h$ 
and $\h/W \to \h/W$ in the commutative diagram in Proposition (3.2). According to the 
notation in the proof of (3.2), we call 
these maps $\tilde{\iota}_{\gamma, \h}$ and  $\iota_{\gamma, \h/W}$ respectively. 
Assume that $S$ is of type $A_r$, $D_r$ or $E_r$.  
When $m = r$, the sheaf $R^2\pi^{an}_*\mathbf{C}$ 
has a trivial monodromy along any $\gamma$.  
In this case, we have $\tilde{\iota}_{\gamma, \h} = id$; hence ${\iota}_{\gamma, \h/W} = id$. 
The problem is when $m < r$. In this case, there is a loop $\gamma$ such that 
$\tilde{\iota}_{\gamma, \h} $ comes from one of  the graph automorphisms 
listed in (4.1). Assume that $\dim \h^ {\tilde{\iota}_{\gamma, \h}} = m$, where 
$\h^{\tilde{\iota}_{\gamma, \h}}$ is the invariant part of $\h$ under $\tilde{\iota}_{\gamma, \h}$.  
By the argument in [Slo, 8.8, Lemma 1], we see that  
$\dim (\h/W)^{\iota_{\gamma, \h/W}} = m$.  Q.E.D.   
\vspace{0.2cm}
 
By using Proposition (4.2), we can prove that the inequality in  
Corollary (1.10) of [Na 1] is actually an equality: 

{\bf Corollary (4.3)}. {\em Let $(X, \omega)$ be a 
projective symplectic variety. Let $U \subset X$ be the locus where 
$X$ is locally a trivial deformation of a 
(surface) rational double point at each $p \in U$. 
Let $\pi: \tilde{U} \to U$ be the minimal resolution and 
let $m$ be the number of irreducible components of 
$\mathrm{Exc}(\pi)$. Then $h^0(U, T^1_U) = m$. } 
\vspace{0.2cm}

{\em Proof} By Lemma (1.5) we obtain a local system $\mathcal{H}$ 
of $\mathbf{C}$-modules as a subsheaf of $T^1_U$. 
Put $\Sigma := \mathrm{Sing}(U)$. Let $\Sigma = \cup \Sigma_i$ be 
the decomposition into connected components. The local system $\mathcal{H}$ 
has support on $\Sigma$. Let $\mathcal{H}_i$ be the restriction 
of $\mathcal{H}$ to each connected component $\Sigma_i$. 
We have an isomorphism: 
$$\mathcal{H} \otimes_{\mathbf{C}}\mathcal{O}_{\Sigma} \cong 
T^1_U.$$  
Then $$h^0(U, T^1_U) = 
h^0(\Sigma, \mathcal{H}\otimes_{\mathbf{C}}\mathcal{O}_{\Sigma}) 
= \Sigma h^0(\mathcal{H}_i)\cdot h^0(\mathcal{O}_{\Sigma_i}).$$ 
Since $\Sigma_i$ can be compactified to a proper normal variety 
$\bar{\Sigma}_i$ such that $\bar{\Sigma}_i - \Sigma_i$ has codimension 
$\geq 2$, we see that $h^0(\mathcal{O}_{\Sigma_i}) = 1$.   Q.E.D. 

\section{Main Results} 

{\bf Theorem (5.1)}. {\em Let $X$ be an affine symplectic 
variety. Then $\mathrm{PD}_X$ is unobstructed.} 
\vspace{0.2cm}

{\em Proof}. (i) Let $U$ be chosen as in (1.8). Let $\pi: \tilde{U} 
\to U$ be the minimal resolution. Put $Z := X \setminus U$. 
In the exact sequence of local cohomology 
$$ ... \to H^i(X, \mathcal{O}_X) \to H^i(U, \mathcal{O}_U) \to 
H^{i+1}_Z(X, \mathcal{O}_X) \to ..., $$ we have $H^{i+1}_Z(X, \mathcal{O}_X) = 
0$ for all $i \le 2$ since $X$ is Cohen-Macaulay and $\mathrm{Codim}_X Z 
\geq 4$. Note that $H^i(X, \mathcal{O}_X) = 0$ for $i > 0$. 
Therefore, one has 
$H^i(U, {\mathcal O}_U) = 0$ for $i = 1,2$. Since $U$ is a symplectic variety, 
$U$ has only rational singularities (cf. (1.1)). In particular, this implies that  
$H^i(\tilde{U}, \mathcal{O}_{\tilde{U}}) = 0$ 
for $i = 1,2$. The resolution $\tilde{U}$ is 
a smooth symplectic variety and $\mathrm{PD}_{\tilde{U}}(\mathbf{C}[\epsilon]) 
\cong H^2(\tilde{U}^{an}, \mathbf{C})$. 
There is a natural map 
$\mathrm{PD}_{\tilde{U}}(\mathbf{C}[\epsilon]) 
\to \mathrm{PD}_U(\mathbf{C}[\epsilon])$. In fact, since 
$R^1\pi_*\mathcal{O}_{\tilde{U}} = 0$ and 
$\pi_*\mathcal{O}_{\tilde{U}} = \mathcal{O}_U$, a first 
order deformation $\tilde{\mathcal{U}}$ 
(without Poisson structure) of $\tilde{U}$ 
induces a first order deformation $\mathcal{U}$ 
of $U$ (cf. [Wa]). Let $\mathcal{U}^0$ be the locus 
where $\mathcal{U} \to \mathrm{Spec}(\mathbf{C}[\epsilon])$ 
is smooth. Since $\tilde{\mathcal{U}} \to \mathcal{U}$ 
is an isomorphism above $\mathcal{U}^0$, the Poisson structure 
of $\tilde{\mathcal{U}}$ induces that of $\mathcal{U}^0$. 
Since the Poisson structure of $\mathcal{U}^0$ uniquely 
extends to that of $\mathcal{U}$, $\mathcal{U}$ becomes a 
Poisson scheme over $\mathrm{Spec}(\mathbf{C}[\epsilon])$. 
This is the desired map. In the same way, one has a 
morphism of functors: 
$$ \mathrm{PD}_{\tilde{U}} \stackrel{\pi_*}\to 
\mathrm{PD}_U. $$ 
Note that $\mathrm{PD}_{\tilde{U}}$ (resp. 
$\mathrm{PD}_U$) has a prorepresentable hull 
$R_{\tilde{U}}$ (resp. $R_U$). Then $\pi_*$ 
induces a local homomorphism of complete local rings: 
$$R_U \to R_{\tilde{U}}.$$ 
We now obtain a 
commutative diagram of exact sequences:  

\begin{equation} 
\begin{CD} 
0 @>>> H^2(U^{an}, \mathbf{C}) @>>> 
\mathrm{PD}_{\tilde{U}}(\mathbf{C}[\epsilon]) @>>> 
H^0(U^{an}, R^2\pi^{an}_*\mathbf{C}) \\ 
@. @V{\cong}VV @VVV @. \\ 
0 @>>> \mathrm{PD}_{lt,U}(\mathbf{C}[\epsilon]) 
@>>> \mathrm{PD}_U(\mathbf{C}[\epsilon]) @>>> 
H^0(\Sigma, \mathcal{H}) 
\end{CD} 
\end{equation} 
  
(ii) Let $E_i$ ($i = 1, ..., m$) be the irreducible components 
of $\mathrm{Exc}(\pi)$. Each $E_i$ defines a class $[E_i] 
\in H^0(U^{an}, R^2\pi^{an}_*\mathbf{C})$. It is easily checked that 
$H^0(U^{an}, R^2\pi^{an}_*\mathbf{C}) = \oplus_{1 \le i \le m}\mathbf{C}[E_i]$.    
This means that 
$$\dim \mathrm{PD}_{\tilde{U}}(\mathbf{C}[\epsilon]) 
= h^2(U^{an}, \mathbf{C}) + m.$$
On the other hand, by Proposition (4.2), 
$h^0(\Sigma, \mathcal{H}) = m$. This means that 
$$\dim \mathrm{PD}_U(\mathbf{C}[\epsilon]) 
\leq h^2(U^{an}, \mathbf{C}) + m.$$ 
As a consequence, we have 
$$ \dim \mathrm{PD}_U(\mathbf{C}[\epsilon]) \leq 
\dim \mathrm{PD}_{\tilde{U}}(\mathbf{C}[\epsilon]).$$
\vspace{0.15cm} 

(iii) We shall prove that the morphism 
$\pi_*: \mathrm{PD}_{\tilde{U}} \to \mathrm{PD}_U$ 
has a finite fiber. More exactly, 
$\mathrm{Spec}(R_{\tilde{U}}) \to 
\mathrm{Spec}(R_U)$ has a finite closed fiber. 
Let $\alpha: R_{\tilde{U}} \to \mathbf{C}[[t]]$ be a 
homomorphism of local $\mathbf{C}$-algebras such 
that the composition map $R_U \to R_{\tilde{U}} 
\stackrel{\alpha}\to \mathbf{C}[[t]]$ is factorized as 
$R_U \to R_U/m_U \to \mathbf{C}[[t]]$. 
Let $U_p$ be the germ of $U^{an}$ at $p \in \Sigma$ 
and let $\tilde{U}_p$ be the germ of 
$\tilde{U}^{an}$ along $(\pi^{an})^{-1}(p)$. 
Denote by $R_{U_p}$ (resp. $R_{\tilde{U}_p}$) 
the prorepresentable hull of the Poisson deformation 
functor $\mathrm{PD}_{U_p}$ (resp. 
$\mathrm{PD}_{\tilde{U}_p}$). 
Since a Poisson deformation of $U$ (resp. $\tilde{U}$) 
induces a Poisson deformation of $U_p$ (resp. $\tilde{U}_p$), 
$\alpha$ induces the map 
$\alpha_p: R_{\tilde{U}_p} \to \mathbf{C}[[t]]$ such that 
the map $R_{U_p} \to R_{\tilde{U}_p} \stackrel{\alpha_p}\to 
\mathbf{C}[[t]]$ is factorizes as $R_{U_p} \to R_{U_p}/m_{U_p} \to 
\mathbf{C}[[t]]$.  
Corresponding to $\alpha$, we have a family of morphisms $\{\pi_n\}_{n \geq 1}$: 
$$ \pi_n: \tilde{U}_n \to U_n, $$ 
where $U_n \cong U \times \mathrm{Spec}\mathbf{C}[t](t^{n+1})$ 
and $\tilde{U}_n$ are Poisson deformations of $\tilde{U}$ 
over $\mathrm{Spec}\mathbf{C}[t]/(t^{n+1})$. 
Restrict these to $\tilde{U}_p$ and $U_p$. Then we have a family 
of morphisms 
$\{\pi_{p,n}\}_{n \geq 1}$: 
$$ \pi_{p,n}: \tilde{U}_{p,n} \to U_{p,n}, $$ 
which are Poisson deformations of $\tilde{U}_p$ and 
$U_p$ determined by $\alpha_p$. 
As proved in (3.1), the map $\mathrm{Spec}(R_{\tilde{U}_p}) 
\to \mathrm{Spec}(R_{U_p})$ is a finite Galois covering. 
This means that each $\tilde{U}_{p,n}$ coincides with the minimal 
resolution of $U_{p,n}$ (i.e. $\tilde{U}_p \times \mathrm{Spec}\mathbf{C}[t]/(t^{n+1})$)  
with the natural Poisson structure determined by that of 
$U_{p,n}$. 
Since all minimal resolution $\tilde{U}_{p,n}$ ($p \in 
\Sigma$) are glued together, we conclude that $\tilde{U}_n \cong  
\tilde{U} \times \mathrm{Spec}\mathbf{C}[t]/(t^{n+1})$ and  
its Poisson structures is 
uniquely determined by that of $U_n$. 
This implies that the given map 
$R_{\tilde{U}} \to \mathbf{C}[[t]]$ factors 
through $R_{\tilde{U}}/m_{\tilde{U}}$. 
\vspace{0.15cm}

(iv) Since the tangent space of $\mathrm{PD}_{\tilde{U}}$ is 
controlled by $H^2(U^{an}, \mathbf{C})$, it has the 
$T^1$-lifting property; hence $\mathrm{PD}_{\tilde{U}}$ 
is unobstructed and $R_{\tilde U}$ is regular. 
\vspace{0.15cm}

(v) By (ii), (iii) and (iv), we conclude that $R_U$ is a regular 
local ring with $\dim R_U = \dim R_{\tilde U}$. 
In fact, since $\dim R_{\tilde U} \le \dim R_U + 
\dim R_{\tilde U}/m_UR_{\tilde U}$, we have $$\dim R_{\tilde U} 
\le \dim R_U$$ by (iii). Since $R_{\tilde U}$ is regular by (iv), 
we have an equality 
$$\dim_{\mathbf C} m_{\tilde U}/(m_{\tilde U})^2 = \dim R_{\tilde U}.$$  
On the other hand, we have an inequality 
$$\dim_{\mathbf C} m_U/(m_U)^2 \geq \dim R_U.$$ 
These three (in)equalities imply that 
$$ \dim_{\mathbf C}m_U/(m_U)^2 \geq \dim_{\mathbf C}m_{\tilde U}/(m_{\tilde U})^2.$$ 
Finally, by (ii), we see that this inequality actually is an equality, and 
the equality 
$\dim R_U = \dim_{\mathbf C}m_U/(m_U)^2$ holds.    
   
Moreover, in the commutative diagram above, 
the map $\mathrm{PD}_U(\mathbf{C}[\epsilon]) \to  
H^0(\Sigma, \mathcal{H})$ is surjective. 
We shall prove that $\mathrm{PD}_X$ is 
unobstructed. Let $S_n := \mathbf{C}[t]/(t^{n+1})$ 
and $S_n[\epsilon] := 
\mathbf{C}[t, \epsilon]/(t^{n+1}, \epsilon^2)$. 
Put $T_n := \mathrm{Spec}(S_n)$ and 
$T_n[\epsilon] := \mathrm{Spec}(S_n[\epsilon])$. 
Let $X_n$ be a Poisson deformation of $X$ over 
$T_n$. Define $\mathrm{PD}(X_n/T_n, T_n[\epsilon])$ 
to be the set of equivalence classes of the 
Poisson deformations of $X_n$ over $T_n[\epsilon]$. 
The $X_n$ induces a Poisson deformation $U_n$ 
of $U$ over $T_n$. Define 
$\mathrm{PD}(U_n/T_n, T_n[\epsilon])$ in a similar 
way. Then, by the same argument as [Na 2, Proposition 13], 
we have $$\mathrm{PD}(X_n/T_n, T_n[\epsilon]) 
\cong \mathrm{PD}(U_n/T_n, T_n[\epsilon]).$$ 
Now, since $\mathrm{PD}_U$ is unobstructed, $\mathrm{PD}_U$ has the 
$T^1$-lifting property. This equality shows that 
$\mathrm{PD}_X$ also has the $T^1$-lifting property. 
Therefore, $\mathrm{PD}_X$ is unobstructed. Q.E.D. 
\vspace{0.2cm} 

(5.2) Let $X$ be an affine symplectic variety. 
Take a (projective) resolution $Z \to X$. By 
Birkar-Cascini-Hacon-McKernan [B-C-H-M], one  
applies the minimal model program to this morphism 
and obtains a relatively minimal model $\pi: Y \to X$.  
The following properties are satisfied:
\vspace{0.2cm}
 
(i) $\pi$ is a crepant, birational projective morphism. 
 
(ii) $Y$ has only $\mathbf{Q}$-factorial terminal 
singularities. 
\vspace{0.2cm}

Note that $Y$ naturally becomes a symplectic variety. 
Let $U \subset X$ be the open locus where, for each 
$p \in U$, the germ $(X,p)$ is non-singular or the 
product of a surface rational double point and a non-singular 
variety. We put $\tilde{U} := \pi^{-1}(U)$.
As in (i) of the proof of Theorem (5.1), the birational maps 
$\pi$ and $\pi\vert_{\tilde U}$ induces natural maps of functors 
$\pi_*: \mathrm{PD}_Y \to \mathrm{PD}_X$ and $(\pi\vert_{\tilde U})_*: 
\mathrm{PD}_{\tilde U} \to \mathrm{PD}_U$. 
There is a commutative diagram of Poisson 
deformation functors 

\begin{equation} 
\begin{CD} 
\mathrm{PD}_Y @>>> \mathrm{PD}_{\tilde U} \\ 
@VVV @VVV \\ 
\mathrm{PD}_X @>>> \mathrm{PD}_U     
\end{CD} 
\end{equation} 
and correspondingly a commutative diagram of prorepresentable  
hulls 

\begin{equation} 
\begin{CD} 
R_{\tilde U} @>>> R_Y  \\ 
@AAA @AAA \\ 
R_U @>>> R_X    
\end{CD} 
\end{equation} 

{\bf Lemma (5.3)}. {\em The horizontal maps $R_{\tilde U} \to R_Y$ 
and $R_U \to R_X$ are both isomorphisms.}  

{\em Proof}. Let $V$ be the regular locus of $Y$. Then $\tilde{U}$ is contained in $V$,  
and we have the restriction map $H^2(V^{an}, \mathbf{C}) \to H^2(\tilde{U}^{an}, \mathbf{C})$.  
This map is an isomorphism by the proof of [Na 3], Proposition 2. 
Note that $\mathrm{PD}_Y(\mathbf{C}[\epsilon]) = H^2(V^{an}, \mathbf{C})$ 
and $\mathrm{PD}_{\tilde{U}}(\mathbf{C}[\epsilon]) = 
H^2(\tilde{U}, \mathbf{C})$. By the $T^1$-lifting principle, 
$\mathrm{PD}_Y$ and $\mathrm{PD}_{\tilde{U}}$ are both unobstructed. 
Let us consider the map $R_{\tilde U} \to R_Y$. By the observation above,   
$R_{\tilde U}$ and $R_Y$ are both regular and the map induces an 
isomorphism of Zariski tangent spaces; hence $R_{\tilde{U}} \cong R_Y$. 
Next let us consider the map $R_U \to R_X$. By Theorem (5.1), 
both local rings are regular and the map induces 
an isomorphism of Zariski tangent spaces; hence $R_U \cong R_X$. Q.E.D. 
\vspace{0.2cm}

By Theorem (5.1), $\dim R_U 
= \dim R_{\tilde{U}}$ and the closed 
fiber of $R_U \to R_{\tilde{U}}$ is finite; hence $\dim R_X = \dim R_Y$ 
and the closed fiber of 
$\pi^*: R_X \to R_Y$ is finite. By the generalized Weierstrass 
preparation theorem, $R_Y$ is a finite $R_X$-module; in other words, 
$\mathrm{Spec}\: R_Y \to \mathrm{Spec}\; R_X$ is a finite morphism.    

We put $R_{X,n} := R_X/m^n$ and $R_{Y,n} := R_Y/(m_Y)^n$. 
Since $\mathrm{PD}_X$ and $\mathrm{PD}_Y$ are both 
prorepresentable, there is a commutative diagram of 
formal universal deformations of $X$ and $Y$: 

\begin{equation} 
\begin{CD} 
\{Y_n\}_{n \geq 1} @>>> \{X_n\}_{n \geq 1}  \\ 
@VVV @VVV \\ 
\mathrm{Spec}(R_{Y,n}) @>>> \mathrm{Spec}(R_{X,n}) 
\end{CD} 
\end{equation} 

(5.4) {\bf Algebraization}   

Let us assume that an affine 
symplectic variety $(X, \omega)$ satisfies the following 
condition (*). 
\vspace{0.2cm}

(*) 

(1)  There is a $\mathbf{C}^*$-action on $X$ 
with only positive weights and a unique fixed point 
$0 \in X$.       

(2) The symplectic form $\omega$ has positive weight $l > 0$.  
\vspace{0.2cm}

By Step 1 of Proposition (A.7) in [Na 2], the $\mathbf{C}^*$-action 
on $X$ uniquely extends to the action on $Y$. 
These $\mathbf{C}^*$-actions induce those on $R_X$ and $R_Y$. 
By Section 4 of [Na 2], $R_Y$ is isomorphic to the 
formal power series ring  $\mathbf{C}[[ y_1, ..., y_d]]$ with 
$wt(y_i) = l$. Since $R_X \subset R_Y$, the $\mathbf{C}^*$-action 
on $R_X$ also has positive weights. 
We put $A := \mathrm{lim} \Gamma (X_n, \mathcal{O}_{X_n})$ and 
$B := \mathrm{lim} \Gamma (Y_n, \mathcal{O}_{Y_n})$. 
Let $\hat{A}$ and $\hat{B}$ be the completions of $A$ and $B$ 
along their maximal ideals. 
Then one has the commutative diagram 

\begin{equation} 
\begin{CD} 
R_X @>>> R_Y  \\ 
@VVV @VVV \\ 
\hat{A} @>>> \hat{B}  
\end{CD} 
\end{equation} 

Let $S$ (resp. $T$) be the $\mathbf{C}$-subalgebra of 
$\hat{A}$ (resp. $\hat{B}$) generated by the eigen-vectors 
of the $\mathbf{C}^*$-action.  On the other hand, 
the $\mathbf{C}$-subalgebra of $R_Y$ generated by eigen-vectors, 
is nothing but $\mathbf{C}[y_1, ..., y_d]$.  Let us consider the 
$\mathbf{C}$-subalgebra of $R_X$ generated by eigen-vectors. 
 By [Na 2], Lemma (A.2), it is generated by eigenvectors 
that form a basis of $m_X/(m_X)^2$.  Since $R_X$ is regular of 
 the same dimension as $R_Y$, the subalgebra is a polynomial 
 ring $\mathbf{C}[x_1, ..., x_d]$.  Now the following commutative diagram 
algebraizes the previous diagram: 

\begin{equation} 
\begin{CD} 
\mathbf{C}[x_1, ..., x_d] @>>> \mathbf{C}[y_1, ..., y_d]  \\ 
@VVV @VVV \\ 
S @>>> T   
\end{CD} 
\end{equation} 
  
By Theorem (5.4.5) of [EGA III], the (formal) birational projective morphism 
$${Y_n} \to \mathrm{Spec}(\hat{B}/(m_{\hat{B}})^n)$$ 
is algebraized to a birational projective morphism 
$$ \hat{Y} \to \mathrm{Spec}(\hat{B}).$$ 
Moreover, by a method similar to that in Appendix of [Na 2], 
this is further algebraized to 
$$ \mathcal{Y} \to \mathrm{Spec}(T). $$   
If we put $\mathcal{X} := \mathrm{Spec}(S)$, then 
we have a $\mathbf{C}^*$-equivariant 
commutative diagram of algebraic schemes 

\begin{equation} 
\begin{CD} 
\mathcal{Y} @>>> \mathcal{X}  \\ 
@VVV @VVV \\ 
\mathrm{Spec}\mathbf{C}[y_1, ..., y_d] @>{\psi}>> \mathrm{Spec}\mathbf{C}[x_1, ..., x_d]   
\end{CD} 
\end{equation} 
\vspace{0.2cm}

{\bf Theorem (5.5)}. {\em In the diagram above, }  
\vspace{0.15cm}

(a) {\em the map $\psi$ is a finite surjective map, }
\vspace{0.15cm}

(b) $\mathcal{Y} \to \mathrm{Spec}\;\mathbf{C}[y_1, ..., y_d]$ {\em is a 
locally trivial deformation of $Y$, and}
\vspace{0.15cm}

(c) {\em the induced birational map $\mathcal{Y}_t \to 
\mathcal{X}_{\psi(t)}$ is an isomorphism for a 
general $t \in \mathrm{Spec}\mathbf{C}[y_1, ..., y_d]$ .} 
\vspace{0.15cm}

{\em Proof}. (a) follows from [Na 2], Lemma (A.4) since $R_Y$ is 
a $R_X$-finite module.    

(b):  Since $Y$ is $\mathbf{Q}$-factorial, $Y^{an}$ is also 
$\mathbf{Q}$-factorial by Proposition (A.9) of [Na 2]. 
Then (b) is Theorem 17 of [Na 2].  

(c) follows from Proposition 24 of [Na 2].     
\vspace{0.2cm}

{\bf  Corollary (5.6)}. {\em Let $(X, \omega)$ be an affine symplectic 
variety with the property (*).   
Then the following two conditions 
are equivalent:} 
\vspace{0.15cm}

(1) {\em $X$ has a crepant projective resolution.}     
\vspace{0.15cm}

(2) {\em $X$ has a smoothing by a Poisson deformation.} 
\vspace{0.2cm}

{\em Proof}. (1) $\Rightarrow$ (2): If $X$ has a crepant resolution, say $Y$. 
By using this $Y$, one can construct a diagram in Theorem (5.5). 
Then, by the property (c), we see that $X$ has a smoothing by 
a Poisson deformation. 

(2) $\Rightarrow$ (1): Let $Y$ be a crepant $\mathbf{Q}$-factorial terminalization 
of $X$. It suffices to prove that $Y$ is smooth.  We again consider 
the diagram in Theorem (5.5). By the assumption, $\mathcal{X}_s$ 
is smooth for a general point $s \in \mathrm{Spec}\mathbf{C}[x_1, ..., x_d]$. 
By the property (a), one can find $t \in   
\mathrm{Spec}\mathbf{C}[y_1, ..., y_d]$ such that $\psi(t) = s$. 
By (c), one has an isomorphism $\mathcal{Y}_t \cong \mathcal{X}_s$. In particular, 
$\mathcal{Y}_t$ is smooth. Then, by (b), $Y (= \mathcal{Y}_0)$ is smooth.

\quad \\
\quad\\

Yoshinori Namikawa \\
Department of Mathematics, 
Faculty of Science, Kyoto University, JAPAN \\
namikawa@math.kyoto-u.ac.jp

\end{document}